\documentclass[aos,MSNbibl,seceqn,nameyear,dvips]{arximspdf}
\usepackage{accents}

\doi{10.1214/13-AOS1110} 
\volume{41}
\issue{3}
\pubyear{2013}
\firstpage{1268}
\lastpage{1298}

\makeatletter

\renewcommand{\mid}{|}

\renewcommand{\underline}{\underaccent{\bar}}

\newcommand{\rright}{\right}
\newcommand{\lleft}{\left}

\newcommand{\cal}{\mathcal}

\newtheorem{theorem}{Theorem}[section]
\newtheorem{lemma}{Lemma}[section]
\newtheorem{cor}{Corollary}[section]

\makeatother

\begin{document}
\begin{frontmatter}

\title{Moment bounds and mean squared prediction errors of long-memory time series\thanksref{T1}}
\runtitle{Moment bounds and prediction errors}

\thankstext{T1}{Supported in part by Grants from HKSAR-RGC-GRF 400410, NSC 100-2118-M-390-001,
NSC 99-2118-M-001-008-MY2 and Academia Sinica Investigator Award.}

\begin{aug}
\author[A]{\fnms{Ngai Hang} \snm{Chan}\ead[label=e1]{nhchan@sta.cuhk.edu.hk}},
\author[B]{\fnms{Shih-Feng} \snm{Huang}\ead[label=e2]{huangsf@nuk.edu.tw}}
\and
\author[C]{\fnms{Ching-Kang} \snm{Ing}\corref{}\ead[label=e3]{cking@stat.sinica.edu.tw}}
\runauthor{N. H. Chan, S.-F. Huang and C.-K. Ing}
\affiliation{Chinese University of Hong Kong,
National University of Kaohsiung,
and Academia Sinica and National Taiwan University}
\address[A]{N. H. Chan\\
Department of Statistics\\
Room 118, Lady Shaw Building\\
Chinese University of Hong Kong\\
Shatin, New Territories\\
Hong Kong\\
\printead{e1}}
\address[B]{S.-F. Huang\\
Department of Applied Mathematics\\
National University of Kaohsiung\\
Kaohsiung, 811\\
Taiwan, ROC\\
\printead{e2}}
\address[C]{C.-K. Ing\\
Institute of Statistical Science\\
Academia Sinica\\
Taipei, 115\\
Taiwan, ROC\\
\printead{e3}} 
\end{aug}

\received{\smonth{12} \syear{2012}}
\revised{\smonth{3} \syear{2013}}

%
\begin{abstract}
A moment bound for the normalized conditional-sum-of-squares (CSS)
estimate of a general autoregressive fractionally integrated moving
average (ARFIMA) model with an arbitrary unknown memory parameter is
derived in this paper. To achieve this goal, a uniform moment bound for
the inverse of the normalized objective function is established. An
important application of these results is to establish asymptotic
expressions for the one-step and multi-step mean squared prediction
errors (MSPE) of the CSS predictor. These asymptotic expressions not
only explicitly demonstrate how the multi-step MSPE of the CSS
predictor manifests with the model complexity and the dependent
structure, but also offer means to compare the performance of the CSS
predictor with the least squares (LS) predictor for integrated
autoregressive models. It turns out that the CSS predictor can gain
substantial advantage over the LS predictor when the integration order
is high. Numerical findings are also conducted to illustrate the
theoretical results.
\end{abstract}

%

\end{frontmatter}

\section{Introduction}\label{sec1}
\label{section1}
Long-memory behavior has been extensively document\-ed in a spectrum of
applications. For background information on
long-memory time series and their applications, readers are referred to
\citet{DouOppTaq03}, 
where important theories and applications of long-memory models in the
areas of finance, insurance, the environment and telecommunications are
surveyed. One distinctive feature of the long-memory phenomenon is that
the autocorrelation function of a long-memory process decays at a
polynomial rate, which is much slower than the exponential rate of a
short-memory process. This feature not only enriches the modeling of
time series data, but also offers new challenges. While considerable
attention has been given in the literature to the derivation of the law
of large numbers and the central limit theorem for the estimated
parameters in many long-memory time series models [see, e.g.,
\citet{Dah89}, \citet{FoxTaq86}, \citet{GirSur90},
\citet{RobHid97} and \citet{Rob06}],
less attention has been devoted to their moment properties.
On the other hand, moment properties of the estimated parameters in
short-memory time series models have been widely studied.
For example, \citet{FulHas81} and \citet{KunYam85}
obtained moment bounds for the least
squares (LS) estimators of stationary autoregressive (AR) models, which
led to asymptotic expressions for the mean
squared prediction error (MSPE) of the corresponding least squares predictors.
\citet{IngWei03} 
established a moment bound for the inverse of Fisher's information
matrix of increasing
dimension under a short-memory AR($\infty$) process, which enabled them
to derive an asymptotic expression for the MSPE of the least
squares predictor of increasing order.
When the moving average (MA) part is taken into account, moment bounds
for the estimated parameters are much more difficult to establish, however.
\citet{ChaIng11} 
recently resolved this difficulty by establishing a \textit{uniform} moment
bound for the inverse of Fisher's information matrix of nonlinear
stochastic regression models.
Based on this bound, they analyzed the MSPE of the
conditional-sum-of-squares (CSS) predictor (defined in Section \ref
{sec3}) and
explained how the final prediction error can be used as an effective
tool in the model selection of autoregressive moving average (ARMA) models.

These aforementioned studies primarily deal with the stationary cases,
which may be inapplicable
in many important situations when nonstationary behaviors are often
encountered.
In view of the importance of incorporating long-memory, short-memory
and nonstationary features simultaneously,
we are led to consider the following general autoregressive
fractionally integrated moving average (ARFIMA) model.
Specifically, suppose the data $y_1,\ldots, y_n$ are generated by
%
%
\begin{eqnarray}
\label{11}
&&
\bigl(1-\alpha_{0,1}B-\cdots-\alpha_{0,p_1}B^{p_1}
\bigr) (1-B)^{d_0}y_t \nonumber\\[-8pt]\\[-8pt]
&&\qquad= \bigl(1-\beta_{0,1}B-
\cdots-\beta_{0,p_2}B^{p_2}\bigr)\varepsilon_{t},\nonumber
\end{eqnarray}
where\vspace*{2pt} $\bolds{\eta}_0=(\bolds{\theta}_0^{\mathsf{T}},d_0)^{\mathsf{T}}=
(\alpha_{0,1},\ldots,\alpha_{0,p_1},\beta_{0,1},\ldots,\beta
_{0,p_2},d_0)^{\mathsf{T}}$ is an unknown coefficient vector with
$d_{0} \in R$ and\vspace*{1pt} $1-\sum_{j=1}^{p_{1}}\alpha_{0,j}z^{j}\neq0$ and
$1-\sum_{j=1}^{p_{2}}\beta_{0,j}z^{j}\neq0$ for $|z|\leq1$, $B$ is the
back-shift operator and $\varepsilon_t$'s are independent random
disturbances with $\mathrm{E}(\varepsilon_t)=0$ for all $t$. Throughout
this paper, it will be assumed that $y_t=\varepsilon_t=0$ for all
$t\le0$. These types of initial conditions are commonly used in the
nonstationary time series literature; see, for example,
\citet{ChaWei88}, \citet{HuaRob11} and \citet{Kat08}.
Assume that
%
%
\begin{equation}
\label{12} \bolds{\theta}_0\times d_0\in\Pi
\times D,
\end{equation}
where
%
%
\begin{equation}
\label{aa} D=[L,U] \qquad\mbox{with } {-}\infty<L<U<\infty
\end{equation}
and
$\Pi$ is a compact set in $R^{p_1+p_2}$ whose element
$\bolds{\theta}=(\alpha_1,\ldots,\alpha_{p_1},\beta_1,\ldots,
\beta_{p_2})^{\mathsf{T}}$ satisfies
%
%
\begin{eqnarray}
\label{13}
&&
A_{1, \bolds{\theta}}(z) = 1-\sum_{j=1}^{p_{1}}
\alpha_{j}z^{j} \neq0,\nonumber\\[-8pt]\\[-8pt]
&&
A_{2, \bolds{\theta}}(z) = 1-\sum
_{j=1}^{p_{2}}\beta_{j}z^{j}
\neq0 \qquad\mbox{for all } |z| \leq1;\nonumber
\\
\label{14}
&&A_{1, \bolds{\theta}}(z) \quad\mbox{and}\quad A_{2, \bolds
{\theta}}(z) \mbox{ have no
common zeros};
\\
\label{15}
&& |\alpha_{p_{1}}|+|\beta_{p_{2}}|>0.
\end{eqnarray}
Note that in the current setting,
$D$ can be any general compact interval in $R$, which
encompasses the important case of nonstationary long-memory models
when \mbox{$d \geq0.5$}.

Let $\varepsilon_t(\bolds{\eta})=A_{1,\bolds{\theta
}}(B)A_{2,\bolds{\theta}}^{-1}(B)(1-B)^d y_t$,
where $\bolds{\eta}=(\eta_{1}, \ldots,\eta_{\bar{p}})^{\mathsf
{T}}=(\bolds{\theta}^{\mathsf{T}}, d)^{\mathsf{T}}$
with $\bar{p}=p_{1}+p_{2}+1$.
Then, the CSS estimate
of $\bolds{\eta}_{0}$, $\hat{\bolds{\eta}}_{n}=(\hat
{\bolds{\theta}}{}^{\mathsf{T}}_{n}, \hat{d}_{n})^{\mathsf{T}}$, is
given by
$\hat{\bolds{\eta}}_{n}=\arg\min_{\bolds{\eta}\in\Pi\times
D}S_n(\bolds{\eta})$,
where $S_n(\bolds{\eta})=\sum_{t=1}^n\varepsilon_t^2(\bolds
{\eta})$
is called the objective\vspace*{1pt} function.
The main goal of this paper is to establish a moment bound for
$n^{1/2}(\hat{\bolds{\eta}}_n-\bolds{\eta}_0)$, namely,
%
%
\begin{equation}
\label{16}
\mathrm{E}\bigl\|n^{1/2}(\hat{\bolds{\eta}}_n-
\bolds{\eta}_0)\bigr\|^q=O(1),\qquad q\ge1,
\end{equation}
where $\| \cdot\|$ denotes the Euclidean norm. We focus on model (\ref{11})
instead of more general ones because of its specific and simple
short-memory component, which makes our proof much more transparent.
On the other hand, it is possible to extend our proof to a broader
class of linear processes;
see the discussion given at the of Section~\ref{sec2} for details.

Although it is assumed in (\ref{11}) that $\mathrm{E}(y_{t})=0$,
this condition is not an issue of overriding concern.
To see this, assume that
$y_{t}=\zeta(t)+A^{-1}_{1, \bolds{\theta
}_{0}}(B)(1-B)^{-d_{0}}A_{2, \bolds{\theta}_{0}}(B)\varepsilon_{t}$
has a mean $\zeta(t)$, where
$\zeta(t)$ is a polynomial in $t$ whose degree $k \in\{0, 1, 2, \ldots
\}$ is known and coefficients are unknown.
Then, it is easy to see that $(1-B)^{k+1} y_{t}$ is a zero-mean ARFIMA
process with memory parameter $d_{0}-k-1$.
Given that (\ref{16}) is valid for any value of $d_{0}$,
the CSS estimate of $\bolds{\eta}^{*}_0=(\bolds{\theta}_{0},
d_{0}-k-1)^{\mathsf{T}}$
based on $(1-B)^{k+1}y_{t}$, say $\hat{\bolds{\eta}}^{*}_n$, still
satisfies
$\mathrm{E}\|n^{1/2}(\hat{\bolds{\eta}}^{*}_n-\bolds{\eta
}^{*}_0)\|^q=O(1), q\ge1$.

An important and interesting consequence of (\ref{16}) is that
asymptotic expressions for the one-step and multi-step
MSPEs of the CSS predictor can be established. These asymptotic
expressions not only explicitly demonstrate how the multi-step MSPE
of the CSS predictor manifests with the model complexity and the
dependent structure, but also offer means to compare the performance
of the CSS predictor with the LS predictor for integrated AR
models.
It is worth mentioning that \citet{HuaRob11} have shown that
$n^{1/2}(\hat{\bolds{\eta}}_n-\bolds{\eta}_0)$ converges in
distribution
to a zero-mean multivariate normal distribution.
However, their result cannot be applied to obtain (\ref{16})
because convergence in distribution does not imply convergence of moments.
While existence of moments of $\hat{\bolds{\eta}}_n$
can be guaranteed easily by the compactness of $\Pi\times D$,
this only yields a bound of $O(n^{q/2})$ for the left-hand side of (\ref{16}),
which is greatly improved by the bound on the right-hand side of (\ref{16}).
Equation (\ref{16}) can also be used to investigate
the higher-order bias and the higher-order efficiency of $\hat
{\bolds{\eta}}_n$.
Because these types of problems require a separate treatment, they are
not pursued in
this paper.

Note that under (\ref{11}) with $d_{0}>1/2$, \citet{Ber95} %
argued that the consistency and asymptotic normality of $\hat
{\bolds{\eta}}_{n}$ should hold.
However, as pointed out by \citet{HuaRob11}, 
the proof given in \citet{Ber95} 
appears to be incomplete
because the property that $\hat{\bolds{\eta}}_{n}$ lies
in a small neighborhood of $\bolds{\eta}_{0}$ with probability
tending to 1
is applied with no justification.
Indeed, this property, reliant on uniform probability bounds for $\{
S_{n}(\bolds{\eta})<S_{n}(\bolds{\eta}_{0})\}$,
is difficult to establish for a general $d$. To circumvent this
difficulty, \citet{HuaRob11} 
partitioned the parameter space (after a small ball
centered at $\bolds{\eta}_{0}$ is removed) into four disjoint subsets
according to the value of $d$, and devised different strategies to
establish uniform probability bounds
for $\{S_{n}(\bolds{\eta})<S_{n}(\bolds{\eta}_{0})\}$
over different subsets. Consequently, the consistency and asymptotic
normality of $\hat{\bolds{\eta}}_{n}$
are first rigorously established in \citet{HuaRob11} %
for model (\ref{11}) with a general~$d$.
However, the uniform probability bounds given in \citet{HuaRob11}, 
converging to zero without rates,
are insufficient to establish (\ref{16}).
To prove (\ref{16}), we would require rates of convergence of uniform
probability bounds,
which are in turn ensured by a uniform moment bound of the inverse of
the normalized objective function,
$a_n^{-1}(d)
\sum_{t=1}^n(\varepsilon_t(\bolds{\eta})-\varepsilon_t(\bolds
{\eta}_0))^2$, where
$a_n(d)=nI_{\{d\ge d_0-1/2\}}+n^{2(d_0-d)}I_{\{d<d_0-1/2\}}$
with $I_{B}$ denoting the indicator function of set $B$.
This uniform moment bound, as stipulated and proved in Lemma \ref
{lemma1}, is based on an argument quite different from those
in \citet{ChaIng11} 
and \citet{HuaRob11}, 
and constitutes one of the major contributions of this article.

In Section~\ref{sec2}, by making use of Lemma~\ref{lemma1} and other uniform
probability/moment bounds, (\ref{16}) is proved in Theorem
\ref{thm21}. The problem of extending (\ref{16}) to a general linear
process that
encompasses (\ref{11}) as a special case is also briefly discussed.\vadjust{\goodbreak} In
Section~\ref{sec3}, Lemma~\ref{lemma1} and Theorem~\ref{thm21} are applied
to derive asymptotic expressions for the one-step and multi-step
MSPEs of the CSS predictors; see Theorems~\ref{thm31} and~\ref{thm33}.
These expressions show that
whereas the contribution of the estimated parameters to the MSPE,
referred to as the second-order MSPE,
in the one-step case only involves the number of the estimated parameters,
the second-order MSPE in the multi-step case reflects more features of
the underlying model,
thereby shedding light about the intriguing
multi-step prediction behaviors of the ARFIMA processes. Another
important implication of Theorems~\ref{thm31} and~\ref{thm33} is that even for an
integrated AR model, the CSS predictor can outperform the LS
predictor when the order of integration is large.
To facilitate the presentation, more technical proofs are deferred to
the \hyperref[app]{Appendix}.
By means of Monte Carlo simulations, we also demonstrate that the finite-sample
behaviors of the one-step and multi-step MSPEs in ARFIMA models can be
revealed by the asymptotic results obtained in Section~\ref{sec3}.
Details of this Monte Carlo study, along with the proof
of (\ref{29}), which is the long-memory counterpart of
Theorem 3.1 of \citet{ChaIng11} and
crucial in proving (\ref{16}),
are provided in the supplementary material [\citet{ChaHuaIng}] in
light of space constraint.

\section{Moment bounds}\label{sec2}
\label{section2}

The major goal of this section is to prove (\ref{16}). To this end, we
need an assumption on $\varepsilon_t$.
\begin{longlist}[(A1)]
\item[(A1)]
There exist $0<\delta_0\le1$, $0<\alpha_0\le1$ and $0<M_1<\infty$ such
that for any $0<s-v\le\delta_0$,
${\sup_{1\le t<\infty, \|{\mathbf v}_t\|=1}}|F_{t,{\mathbf
v}_t}(s)-F_{t,{\mathbf
v}_t}(v)|\le M_1(s-v)^{\alpha_0}$,
where ${\mathbf v}_t$ is a $t$-dimensional vector and $F_{t,{\mathbf
v}_t}(\cdot
)$ denotes the distribution of
${\mathbf v}_t^{\mathsf{T}}(\varepsilon_t,\ldots,\varepsilon_1)$.
\end{longlist}
Note that an assumption like (A1) has been used in the literature to
deal with
the moment properties of the LS estimates in the AR or ARMA context;
see, for example,
\citet{FinWei02}, \citet{Ing03}, \citet{Sch05} and \citet{ChaIng11}.
When $\varepsilon_{t}$'s are normally distributed,
(A1) is satisfied with $M_{1}=(2 \pi\sigma^{2})^{-1/2}$ and $\alpha
_{0}=1$ for any $\delta_{0}>0$.
In addition, when $\varepsilon_{t}$'s are i.i.d. with an integrable
characteristic function,
(A1) is satisfied with any $\delta_{0}>0$, $\alpha_{0}=1$ and some $M_{1}>0$.
For a more detailed discussion of (A1), see \citet{IngSin06}. %

The following two lemmas, which may be of independent interest, play a
key role in proving (\ref{16}).
Let $B_{\delta}(\bolds{\eta}_0)=\{\bolds{\eta} \in R^{\bar
{p}}\dvtx \|\bolds{\eta}-\bolds{\eta}_{0}\|<\delta\}$.

%
\begin{lemma}
\label{lemma1}
Assume (\ref{11})--(\ref{15}) and \textup{(A1)}. Then, for any $\delta>0$ such
that $\Pi\times D-B_{\delta}(\bolds{\eta}_0)\ne\varnothing$,
any $q>0$ and any $\theta>0$, we have
%
%
\begin{equation}
\label{21}\quad \mathrm{E} \Biggl[ \Biggl\{\inf_{\bolds{\eta}\in\Pi\times
D-B_{\delta}(\bolds{\eta}_0)}a_n^{-1}(d)
\sum_{t=1}^n\bigl(\varepsilon_t(
\bolds{\eta})-\varepsilon_t(\bolds{\eta}_0)
\bigr)^2 \Biggr\}^{-q} \Biggr] =O \bigl((\log
n)^{\theta} \bigr).
\end{equation}
\end{lemma}
To perceive\vspace*{1pt} the subtlety of Lemma~\ref{lemma1}, first express
$a_n^{-1}(d) \sum_{t=1}^n(\varepsilon_t(\bolds{\eta})-\varepsilon
_t(\bolds{\eta}_0))^2$
as $n^{-1} \sum_{t=1}^n g^{2}_t(\bolds{\eta})$, where\vadjust{\goodbreak}
$g_t(\bolds{\eta})=g_{t,n}(\bolds{\eta})=n^{1/2}
(\varepsilon_t(\bolds{\eta})-\varepsilon_t(\bolds{\eta}_0))
\times a_n^{-1/2}(d)$.
Since $g_{t}(\bolds{\eta})$ is a scalar-valued continuous function
on $\Pi\times D-B_{\delta}(\bolds{\eta}_0)$,
in view of the proof of Theorem 2.1 of \citet{ChaIng11},
(\ref{21})~follows if we can show that
$g_{t}(\bolds{\eta})$ satisfy conditions (C2) and (C3) of the same paper
with slight modifications to accommodate the triangular array feature
of $g_{t}(\bolds{\eta})$.
However, for $d \leq d_{0}-1/2$ and for all large $n$,
the correlation between $g_{t}(\bolds{\eta})$ and
$g_{s}(\bolds{\eta})$
is overwhelmingly large if $t, s \to\infty$ as $n \to\infty$ and $|t-s|$
is bounded by a positive constant.
Therefore, even when (A1) is imposed,
it is still difficult to find
a positive constant $b$ such that for all large $t$ and $n$,
the conditional distribution of $g_{t}(\bolds{\eta})$
given $\{\varepsilon_{s}, s \leq t-b\}$ is sufficiently smooth,
which corresponds to (C2) of \citet{ChaIng11}. 
Moreover,
while $g_{t}(\bolds{\eta})$ is continuous on $\Pi\times
D-B_{\delta}(\bolds{\eta}_0)$,
it is not differentiable at $d=d_{0}-1/2$,
making it quite cumbersome to prove
that there exist $c_{1}>0$
and nonnegative random variables $B_{t}$'s satisfying
$\max_{1\leq t \leq n} \mathrm{E}(B_{t}) =O(1)$ such that for \textit{all}
$\bolds{\eta}_{1}, \bolds{\eta}_{2} \in\Pi\times D-B_{\delta
}(\bolds{\eta}_0)$
with $\|\bolds{\eta}_{1}-\bolds{\eta}_{2}\|<c_{1}$,
$
|g_{t}(\bolds{\eta}_{1})-g_{t}(\bolds{\eta}_{2})| \leq B_{t}\|
\bolds{\eta}_{1}-\bolds{\eta}_{2}\|$ a.s.,
which corresponds to (C3) of \citet{ChaIng11}.
Indeed, this latter condition is particularly difficult to verify
when $\bolds{\eta}_{1}$ and $\bolds{\eta}_{2}$ lie on
\textit{different sides} of the hyperplane $d=d_{0}-1/2$.
As will be seen in the \hyperref[app]{Appendix}, the
$B_{t}$'s derive in (\ref{A7}) and (\ref{A8})
no longer satisfy $\max_{1\leq t \leq n} \mathrm{E}(B_{t}) =O(1)$,
which also results in a slowly varying component
on the right-hand side of~(\ref{21}).

Throughout this paper, $C$ represents a generic positive constant,
independent of~$n$,
whose value may differ from one occurrence to another.

%
\begin{lemma}
\label{lemma2}
Assume (\ref{11})--(\ref{15}), \textup{(A1)} and
%
%
\begin{equation}
\label{22} \sup_{t\ge1}\mathrm{E}|\varepsilon_t|^{q_1}<
\infty,
\end{equation}
where $q_1>q\ge2$.
Let $\delta$ satisfy $\Pi\times D-B_{\delta}(\bolds{\eta}_0)\ne
\varnothing$ and $v>0$ be a small constant.
Define
$B_{0,v}=\{(\bolds{\theta}^{\mathsf{T}},d)^{\mathsf{T}}\dvtx
(\bolds{\theta}^{\mathsf{T}},d)^{\mathsf{T}}\in\Pi\times
D-B_{\delta}(\bolds{\eta}_0)$
with $d_0-\frac12\le d\le U\}$
and
$B_{j,v}=\{(\bolds{\theta}^{\mathsf{T}},d)^{\mathsf{T}}\dvtx
(\bolds{\theta}^{\mathsf{T}},d)^{\mathsf{T}}\in\Pi\times
D-B_{\delta}(\bolds{\eta}_0)$
with $d_0-1/2-jv\le d\le d_0-1/2-(j-1)v\}$, $j\ge1$.
Then, for any $\theta>0$,
%
%
\begin{eqnarray}
\label{23}
&&
\mathrm{E} \biggl\{\frac{{\sup_{\bolds{\eta}\in B_{j,v}}}
|{\sum_{t=1}^n}(\varepsilon_t(\bolds{\eta})-\varepsilon
_{t}(\bolds{\eta}_0))\varepsilon_t| } {
\inf_{\bolds{\eta}\in B_{j,v}}\sum_{t=1}^n
(\varepsilon_t(\bolds{\eta})-\varepsilon_t(\bolds{\eta
}_0))^2} \biggr\}^q \nonumber\\[-8pt]\\[-8pt]
&&\qquad \leq \cases{\displaystyle C \biggl(\frac{(\log n)^{3/2}}{n^{1/2}}
\biggr)^{q}(\log n)^\theta, &\quad$j=0$, \vspace*{2pt}\cr
\displaystyle C \biggl(\frac{\log
n}{n^{1/2+jv-2v}} \biggr)^q(\log n)^\theta, &\quad$j\ge1$.}\nonumber
\end{eqnarray}
\end{lemma}

Lemma~\ref{lemma2} implies
\begin{eqnarray*}
P(\hat{\bolds{\eta}}_{n} \in B_{j,v})&\leq&
P\Bigl(\inf_{\bolds{\eta}\in B_{j,v} }S_{n}(\bolds{\eta})\leq
S_{n}(\bolds{\eta}_{0})\Bigr)
=O\bigl(\bigl\{(\log n)^{3/2}/n^{1/2}\bigr\}^{q}(\log n)^\theta\bigr)
\end{eqnarray*}
for $j=0$,
and $O(\{\log n/n^{1/2+jv-2v}\}^{q}(\log n)^\theta)$ for $j\geq1$,
suggesting that for $d<d_{0}-1/2$, the smaller the value of $d$, the
less likely $\hat{d}_{n}$ will fall in a neighborhood of~$d$. These
probability bounds can suppress the orders of magnitude of
$\|n^{1/2}(\hat{\bolds{\eta}}_n-\bolds{\eta}_0)\|^q$ and
$\sup_{\bolds{\eta} \in
B_{j,v}}n\{\varepsilon_{n+1}(\bolds{\eta})-\varepsilon
_{n+1}(\bolds{\eta}_0)\}^2$,
thereby\vspace*{-2pt} yielding that $\mathrm{E}\{\|n^{1/2}(\hat{\bolds{\eta}}_n-\bolds
{\eta}_0)\|^q \times
I_{\{\hat{\bolds{\eta}}_n\in\Pi\times D-B_{\delta}(\bolds{\eta
}_0)\}} \}$
and\vspace*{1pt} $n\mathrm{E}[\{\varepsilon_{n+1}(\hat{\bolds{\eta}}_n)-\varepsilon
_{n+1}(\bolds{\eta}_0)\}^2
I_{\{\hat{\bolds{\eta}}_n\in\Pi\times D-B_{\delta}(\bolds{\eta
}_0)\}}]$
are asymptotically negligible; see Corollary~\ref{cor21} and Lemma \ref
{lemma31}.
As will become clear later, the first moment property is indispensable
for proving (\ref{16}),
whereas the second one is important in analyzing the MSPE of the CSS predictor.
It is also worth mentioning that
the order of magnitude of
$\sup_{\bolds{\eta} \in B_{j,v}}n\{\varepsilon_{n+1}(\bolds
{\eta})-\varepsilon_{n+1}(\bolds{\eta}_0)\}^2$
is
$n (\log n)^{3}$ for $j=0$ and $n^{1+2vj} (\log n)^{2}$ for $j\geq1$,
which increases as $j$ does; see (\ref{35}) for more details.
The next corollary is a direct application of Lemma~\ref{lemma2}.

%
\begin{cor}
\label{cor21}
Suppose that the assumptions of Lemma~\ref{lemma2} hold. Then, for any
$\delta>0$ such that
$\Pi\times D-B_{\delta}(\bolds{\eta}_0)\ne\varnothing$,
%
%
\begin{equation}
\label{24} \mathrm{E} \bigl\{\bigl\|n^{1/2}(\hat{\bolds{
\eta}}_n-\bolds{\eta}_0)\bigr\|^q
I_{\{\hat{\bolds{\eta}}_n\in\Pi\times D-B_{\delta}(\bolds{\eta
}_0)\}} \bigr\}=o(1).
\end{equation}
\end{cor}
\begin{pf}
Since both $\hat{\bolds{\eta}}_n$ and $\bolds{\eta}_0$ are in
$\Pi\times D$,
$\|\hat{\bolds{\eta}}_n\|$ and $\|\bolds{\eta}_0\|$ are
bounded above by a finite constant.
Therefore, it suffices for (\ref{24}) to show that
%
%
\begin{equation}
\label{25} P \bigl(\hat{\bolds{\eta}}_n\in\Pi\times
D-B_{\delta}(\bolds{\eta}_0) \bigr)=o
\bigl(n^{-q/2}\bigr).
\end{equation}
Let $q_1>q_1^*>q$ and $0<v<\frac12(1-q/q_1^*)$. Without loss of
generality, assume that $L=d_0-(1/2)-Wv$ for some large integer $W>0$.
Then, it follows from Lem\-ma~\ref{lemma2} (with $q=q_1^*$) and
Chebyshev's inequality that for any $\theta>0$,
\begin{eqnarray*}
P \bigl(\hat{\bolds{\eta}}_n\in\Pi\times D-B_{\delta}(
\bolds{\eta}_0) \bigr) &\leq&\sum_{j=0}^W
P(\hat{\bolds{\eta}}_n\in B_{j,v})
\\
&\leq& C \sum_{j=0}^W \mathrm{E} \biggl\{
\frac{{\sup_{\bolds{\eta}\in B_{j,v}}}
|{\sum_{t=1}^n}(\varepsilon_t(\bolds{\eta})-\varepsilon
_{t}(\bolds{\eta}_0))\varepsilon_t| } {
\inf_{\bolds{\eta}\in B_{j,v}}\sum_{t=1}^n
(\varepsilon_t(\bolds{\eta})-\varepsilon_t(\bolds{\eta
}_0))^2} \biggr\}^{q_{1}^{*}}
\\
&\leq& C \biggl(\frac{\log n}{n^{1/2-v}} \biggr)^{q_1^*}(\log
n)^{\theta}\\
&=&o\bigl(n^{-q/2}\bigr),
\end{eqnarray*}
which gives (\ref{25}).
\end{pf}

While Theorem 2.1 of \citet{HuaRob11} 
showed that $\hat{\bolds{\eta}}_n \to_{p} \bolds{\eta}_0$
under substantially weaker assumptions on $\varepsilon_{t}$, it seems
tricky to extend their arguments
to obtain a convergence rate like the one given in (\ref{25}), which
is critical to proving of (\ref{24}).
As a by-product of (\ref{25}), we obtain
$\hat{\bolds{\eta}}_n \to\bolds{\eta}_0$ a.s., which
follows immediately from (\ref{25}) with $q>2$
and the Borel--Cantelli lemma.\vadjust{\goodbreak}
The main result is given in the next theorem.
First, some notation. For $1\leq m\leq\bar{p}$,
define $\mathbf{J}(m,\bar{p})=\{(j_{1},\ldots, j_{m})\dvtx
j_{1}<\cdots<j_{m}, j_{i}\in\{1,\ldots, \bar{p}\}, 1\leq i \leq m \}
$, and for $\mathbf{j}=(j_{1},\ldots, j_{m})
\in\mathbf{J}(m,\bar{p})$ and a smooth function
$w=w(\xi_{1}, \ldots, \xi_{\bar{p}})$, let $\mathbf{D}_{\mathbf{j}} w=
\partial^{m} w/\partial\xi_{j_{1}},\ldots,
\partial\xi_{j_{m}}$.

%
\begin{theorem}
\label{thm21}
Assume (\ref{11})--(\ref{15}), \textup{(A1)},
%
%
\begin{equation}
\label{26} \sup_{t \geq1}{\mathrm E}|\varepsilon_{t}|^{4q_{1}}<
\infty,\qquad q_{1}>q\geq1,
\end{equation}
and
%
%
\begin{equation}
\label{27} \bolds{\eta}_{0} \in\operatorname{int} \Pi\times D.
\end{equation}
Then (\ref{16}) holds.
\end{theorem}

\begin{pf} Since by (\ref{25}) or Theorem 2.1 of \citet{HuaRob11}, 
\mbox{$\hat{\bolds{\eta}}_{n} \to\bolds{\eta}_{0}$} in probability, and
since (\ref{27}) is assumed, there exists $0<\tau_{1}<\{1-(q/q_{1})\}
/2$ such
that
%
%
\begin{equation}
\label{28} B_{\tau_{1}}(\bolds{\eta}_{0}) \subset\Pi
\times D \quad\mbox{and}\quad\lim_{n \to\infty} P\bigl(\hat{\bolds{
\eta}}_{n} \in B_{\tau_{1}}(\bolds{\eta}_{0})
\bigr)=1.
\end{equation}
Let
$\nabla\varepsilon_t(\bolds{\eta})=
(\nabla\varepsilon_t(\bolds{\eta})_1,\ldots,
\nabla\varepsilon_t(\bolds{\eta})_{\bar{p}})^{\mathsf{T}}
=(\partial\varepsilon_t(\bolds{\eta})/\partial\eta_1,\ldots,
\partial\varepsilon_t(\bolds{\eta})/\partial\eta_{\bar
{p}})^{\mathsf{T}}$
and $\nabla^2\varepsilon_t(\bolds{\eta})=
(\nabla^2\varepsilon_t(\bolds{\eta})_{i,j})=
(\partial^{2} \varepsilon_t(\bolds{\eta})/\partial\eta
_{i}\,\partial\eta_{j})$.
Assume first that the following relations hold:
%
%
\begin{eqnarray}
\label{29} \mathrm{E}\Biggl\{\sup_{\bolds{\eta}\in B_{\tau
_1}(\bolds{\eta}_0)} \lambda_{\min}^{-\gamma}
\Biggl(n^{-1}\sum_{t=1}^n \nabla
\varepsilon_t(\bolds{\eta}) \bigl(\nabla\varepsilon_t(
\bolds{\eta})\bigr)^{\mathsf{T}}\Biggr)\Biggr\}&=&O(1)\nonumber\\[-9pt]\\[-9pt]
&&\eqntext{\mbox{for
any } \gamma
\geq1,}
\\[-2pt]
%
\label{210} \max_{1\le i,j\le\bar{p}}\mathrm{E}\Biggl\{\sup
_{\bolds{\eta}\in B_{\tau_1}(\bolds{\eta}_0)} \Biggl|n^{-1/2}\sum_{t=2}^n
\varepsilon_t\nabla^2\varepsilon_t(
\bolds{\eta})_{i,j}\Biggr|^{q_1}\Biggr\}&=&O(1),
\\[-2pt]
%
\label{211} \max_{1\le i,j\le\bar{p},2\le t\le n}\mathrm{E}\Bigl
\{\sup
_{\bolds{\eta}\in B_{\tau_1}(\bolds{\eta}_0)} \bigl|\nabla
^2\varepsilon_t(
\bolds{\eta})_{i,j}\bigr|^{4q_1}\Bigr\}&=&O(1),
\\[-2pt]
%
\label{212} \max_{1\le i\le\bar{p}, 2\le t\le n}\mathrm{E}\Bigl\{
\sup
_{\bolds{\eta}\in B_{\tau_1}(\bolds{\eta}_0)} \bigl|\nabla
\varepsilon_t(\bolds{
\eta})_{i}\bigr|^{4q_1}\Bigr\}&=&O(1),
\\[-2pt]
%
\label{213} P\Biggl\{\sup_{\bolds{\eta}\in B_{\tau_1}(\bolds
{\eta}_0)} \lambda_{\min}^{-1}
\Biggl(n^{-1}\sum_{t=2}^n\nabla
\varepsilon_t(\bolds{\eta}) \bigl(\nabla\varepsilon_t(
\bolds{\eta})\bigr)^{\mathsf{T}}\Biggr) >\bar{M}\Biggr\}&=&O
\bigl(n^{-q}\bigr)\nonumber\\[-9pt]\\[-9pt]
&&\eqntext{\mbox{for some } \bar{M}>0,}
\\[-2pt]
%
\label{214} P\Biggl\{\sup_{\bolds{\eta}\in B_{\tau_1}(\bolds
{\eta}_0)} n^{-1}\sum
_{t=2}^n \bigl\|\nabla\varepsilon_t(
\bolds{\eta})\bigr\|^2 >\bar{M}\Biggr\}&=&O\bigl(n^{-q}\bigr)\nonumber\\[-9pt]\\[-9pt]
&&\eqntext{\mbox{for some } \bar{M}>0,}
\\
%
\label{215} \max_{1\le i,j\le\bar{p}} P\Biggl\{\sup_{\bolds{\eta
}\in B_{\tau_1}(\bolds{\eta}_0)}
n^{-1}\sum_{t=2}^n \bigl(
\nabla^2\varepsilon_t(\bolds{\eta})_{i,j}
\bigr)^2 >\bar{M}\Biggr\}&=&O\bigl(n^{-q}\bigr)\nonumber\\[-8pt]\\[-8pt]
&&\eqntext{\mbox{for
some } \bar{M}>0.}
\end{eqnarray}
Then, making use of (\ref{28})--(\ref{215}) and an argument given in
the proof of Theorem~2.2 of \citet{ChaIng11}, 
we obtain
%
%
\begin{equation}
\label{216} \mathrm{E}\bigl(\bigl\|n^{1/2}(\hat{\bolds{
\eta}}_{n} - \bolds{\eta}_{0})\bigr\|^{q}I_{O_{1n}}
\bigr) = O(1),
\end{equation}
where
$O_{1n}=\{\hat{\bolds{\eta}}_{n} \in B_{\tau_{1}^{*}}
(\bolds{\eta}_{0})\}$ with
$0<\tau^{*}_{1} < \min\{\tau_{1}, 3^{-1}\bar{p}^{-1}\bar{M}^{-2}\}$.
Moreover, it follows from Corollary~\ref{cor21} that
%
%
\begin{equation}
\label{217} \mathrm{E}\bigl(\bigl\|n^{1/2}(\hat{\bolds{
\eta}}_{n} - \bolds{\eta}_{0})\bigr\|^{q}I_{O_{2n}}
\bigr) = O(1),
\end{equation}
where $O_{2n}=\{\hat{\bolds{\eta}}_{n} \in\Pi\times D-B_{\tau_{1}^{*}}
(\bolds{\eta}_{0})\}$.
Combining (\ref{216}) and (\ref{217}) gives the desired conclusion
(\ref{16}).
To complete the proof, it remains to show that (\ref{29})--(\ref
{215}) are true.

A proof of (\ref{29}), which is similar to that of Theorem 3.1 of \citet{ChaIng11}, 
but needs to be modified with the long-memory effect of
$\nabla\varepsilon_{t}(\bolds{\eta})$, $\bolds{\eta} \in
B_{\tau_{1}}(\bolds{\eta}_{0})$,
is deferred to the supplementary material [\citet{ChaHuaIng}].
To prove (\ref{215}), write
\begin{eqnarray*}
\varepsilon_t(\bolds{\eta})&=&(1-B)^{d-d_0}A_{2,\bolds{\theta
}_0}(B)A_{1,\bolds{\theta}_0}^{-1}(B)
A_{1,\bolds{\theta}}(B)A_{2,\bolds{\theta}}^{-1}(B)\varepsilon
_t\\
&=&\sum_{s=0}^{t-1}b_{s}(\bolds{\eta})\varepsilon_{t-s},
\end{eqnarray*}
where $b_{0}(\bolds{\eta})=1$.
Then, with
$c_{s,ij}(\bolds\eta)= \partial^{2} b_{s}(\bolds\theta
)/\partial\eta_{i}\,\partial\eta_{j}$
and $b_{s,i}(\bolds\eta)= \partial b_{s}(\bolds\eta)/\break\partial
\eta_{i}$,
$\nabla\varepsilon_t(\bolds{\eta})_{i}
=\sum_{s=1}^{t-1}b_{s,i}(\bolds\eta)\varepsilon_{t-s}$ and
$\nabla^2\varepsilon_t(\bolds{\eta})_{i,j}=
\sum_{s=1}^{t-1}c_{s,ij}(\bolds{\eta})\varepsilon_{t-s}$.
It is clear that $b_{s,i}(\bolds\eta)$ and $c_{s,ij}(\bolds
\eta)$ have continuous partial derivatives,
$\mathbf{D}_{\mathbf{j}}b_{s,i}(\bolds\eta)$ and $\mathbf
{D}_{\mathbf{j}}c_{s,ij}(\bolds\eta)$,
on $B_{\tau_{1}}(\bolds\eta_{0})$. Moreover, it follows from
arguments similar to those in the proofs of
Theorem 4.1 of \citet{Lin07} 
and Lemma 4 of \citet{HuaRob11} 
that for any $s\geq1$,
%
%
\begin{equation}
\label{218} \max_{1\leq i,j\leq\bar{p}}\sup_{\bolds{\eta} \in
B_{\tau_{1}(\bolds{\eta}_{0})}}\bigl|c_{s, ij}(
\bolds{\eta})\bigr| \leq C \bigl(\log(s+1)\bigr)^{2} s^{-1+\tau_{1}}
\end{equation}
and
%
%
\begin{equation}
\label{219}\qquad \max_{1\leq i,j\leq\bar{p}}\max_{\mathbf{j} \in\mathbf
{J}(m, \bar{p}), 1\leq m \leq\bar{p}}\sup
_{\bolds{\eta} \in B_{\tau_{1}(\bolds{\eta}_{0})}}\bigl|D_{\mathbf
{j}}c_{s, ij} (\bolds{\eta})\bigr|
\leq C \bigl(\log(s+1)\bigr)^{3} s^{-1+\tau_{1}}.
\end{equation}
Define
\begin{eqnarray*}
A_{s}(i,j)&=&\max_{\mathbf{j} \in\mathbf{J}(m,\bar{p}),1\leq m \leq\bar
{p} }\sup_{\bolds{\eta}\in
B_{\tau_{1}}(\bolds{\eta}_{0})}\bigl(\mathbf{D}_{\mathbf
{j}}c_{s,ij}(\bolds{\eta})\bigr)^{2},\\
S_{r,s}(i,j)&=&
\max_{{\mathbf j}\in{\mathbf J}(m, \bar{p}), 1\leq m \leq\bar{p}}\sup
_{\bolds{\eta} \in B_{\tau_{1}}(\bolds{\eta}_{0})}
\bigl|D_{\mathbf{j}}\bigl\{c_{r, ij}(\bolds{\eta})c_{s, ij}(\bolds{\eta})\bigr\}\bigr|
\end{eqnarray*}
and $B_{s}(i,j)= \sup_{\bolds{\eta}\in
B_{\tau_{1}}(\bolds{\eta}_{0})}c^{2}_{s,ij}(\bolds{\eta})$. Equations
(\ref{218}) and (\ref{219}) yield that for any $1\leq i, j\leq\bar{p}$,
%
%
\begin{eqnarray}
\label{220} \sum_{l=1}^{\infty}c^{2}_{l, ij}(
\bolds{\eta}_{0}) &\leq& C,\nonumber\\[-8pt]\\[-8pt]
\Biggl\{\sum
_{l=1}^{\infty}S_{l,l}(i,j)\Biggr
\}^{2} &\leq& C\Biggl[ \Biggl\{\sum_{l=1}^{\infty}
A_{l}(i,j)\Biggr\}^{2} +\Biggl\{\sum
_{l=1}^{\infty}B_{l}(i,j)\Biggr
\}^{2}\Biggr] \leq C.\nonumber
\end{eqnarray}
On the other hand, by (B.6) of \citet{ChaIng11}, 
Chebyshev's inequality and (\ref{26}),
we have for
$\bar{M}>2\sigma^2 \max_{1\leq i, j \leq r}\sum_{s=1}^\infty\sup
_{\bolds{\eta}\in B_{\tau_{1}}(\bolds{\eta}_{0})}
c^2_{s,ij}(\bolds{\eta})$ and any $1\leq i$, $j \leq\bar{p}$,
%
%
\begin{eqnarray}
\label{221}
&&
P\Biggl\{\sup_{\bolds{\eta}\in B_{\tau
_1}(\bolds{\eta}_0)} n^{-1}\sum
_{t=2}^n\bigl(\nabla^2
\varepsilon_t(\bolds{\eta})_{i,j}\bigr)^2
>\bar{M}\Biggr\}
\nonumber\hspace*{-35pt}
\\
&&\quad\le P\Biggl\{ \sup_{\bolds{\eta}\in B_{\tau_1}(\bolds{\eta
}_0)} \Biggl| n^{-1}\sum
_{t=2}^n\bigl[ \bigl(\nabla^2
\varepsilon_t(\bolds{\eta})_{i,j}
\bigr)^2- \mathrm{E}\bigl(\nabla^2\varepsilon_t(
\bolds{\eta})_{i,j}\bigr)^2 \bigr]\Biggr|^{2q_{1}}\!\! >\!
\biggl(\frac{\bar{M}}{2}\biggr)^{2q_{1}}\! \Biggr\}
\nonumber\hspace*{-35pt}\\
&&\quad\le Cn^{-2q}\Biggl\{ \sum_{s=1}^{n-1}
\Biggl(\sum_{l=1}^{n-s}c^{2}_{l, ij}(
\bolds{\eta}_{0})\Biggr)^{2}+ \sum
_{s=1}^{n-1}\Biggl(\sum_{l=1}^{n-s}S_{l,l}(i,j)
\Biggr)^{2}\Biggr\}^{q_{1}}\hspace*{-35pt}
\\
&&\qquad{}+ Cn^{-q_{1}-1} \sum_{r=2}^{n-1}
\Biggl\{ \Biggl[\sum_{s=1}^{r-1}\Biggl(\sum
_{l=1}^{n-s}S_{l+r-s, l}(i,j)
\Biggr)^{2}\Biggr]^{q_{1}}
\nonumber\hspace*{-35pt}
\\
&&\hspace*{76pt}\qquad{}+ \Biggl[\sum_{s=1}^{r-1}\Biggl(\sum
_{l=1}^{n-s}\bigl|c_{l+r-s, ij}(\bolds{
\eta}_{0}) c_{l, ij}(\bolds{\eta}_{0})\bigr|
\Biggr)^{2}\Biggr]^{q_{1}}\Biggr\}.
\nonumber\hspace*{-35pt}
\end{eqnarray}
Since (\ref{218}) and (\ref{219}) also ensure that for some\vspace*{2pt} $\tau
_{1}< \tau_{2}< \{1-(q/q_{1})\}/2$, any $1\leq i, j \leq\bar{p}$
and any $r>s$,
$\{\sum_{l=1}^{\infty}S_{l+r-s, l}(i,j)\}^{2} \leq
C(r-s)^{-1+2\tau_{2}}$\vspace*{2pt}
and
$(\sum_{l=1}^{\infty}|c_{l+r-s, ij}(\bolds{\eta}_{0})c_{l,
ij}(\bolds{\eta}_{0})|)^{2}
\leq C(r-s)^{-1+2\tau_{2}}$,\vspace*{2pt}
we conclude from these, (\ref{220}) and (\ref{221}) that for any
$1\leq i, j\leq\bar{p}$,
\[
P\Biggl( \sup_{\bolds{\eta}\in B_{\tau_1}(\bolds{\eta}_0)}
n^{-1}\sum
_{t=2}^n\bigl(\nabla^2
\varepsilon_t(\bolds{\eta})_{i,j}\bigr)^2
>\bar{M}\Biggr)=O\bigl(n^{-q_{1}(1-2\tau_{2})}\bigr)=o\bigl(n^{-q}\bigr).
\]
Thus, (\ref{215}) follows.

By analogy with (\ref{218}) and (\ref{219}), we have
%
%
\begin{equation}
\label{222} \max_{1\leq i \leq\bar{p}}\sup_{\bolds{\eta} \in
B_{\tau_{1}(\bolds{\eta}_{0})}}\bigl|b_{s, i}(
\bolds{\eta})\bigr| \leq C \bigl(\log(s+1)\bigr) s^{-1+\tau_{1}}
\end{equation}
and
%
%
\begin{equation}
\label{223}\qquad \max_{1\leq i\leq\bar{p}}\max_{\mathbf{j} \in\mathbf
{J}(m, \bar{p}), 1\leq m \leq\bar{p}}\sup
_{\bolds{\eta} \in B_{\tau_{1}(\bolds{\eta}_{0})}}\bigl|D_{\mathbf
{j}}b_{s, i} (\bolds{\eta})\bigr|
\leq C \bigl(\log(s+1)\bigr)^{2} s^{-1+\tau_{1}}.\vadjust{\goodbreak}
\end{equation}
In addition, (0.3) and (0.4) in the supplementary material
[\citet{ChaHuaIng}] ensure that there exists $\underline{c}>0$
such that for all large $n$,
%
%
\begin{equation}
\label{224} \inf_{\bolds{\eta}\in B_{\tau_1}(\bolds{\eta
}_0)} \lambda_{\min}\Biggl(
n^{-1}\sum_{t=2}^n \mathrm{E}
\bigl\{\nabla\varepsilon_t(\bolds{\eta}) \bigl(\nabla
\varepsilon_t(\bolds{\eta})\bigr)^{\mathsf{T}}\bigr\}\Biggr)>
\underline{c}.
\end{equation}
Denote $\nabla\varepsilon_t(\bolds{\eta})(\nabla\varepsilon
_t(\bolds{\eta}))^{\mathsf{T}}$
and $ \mathrm{E}\{\nabla\varepsilon_t(\bolds{\eta})(\nabla
\varepsilon_t(\bolds{\eta}))^{\mathsf{T}}\}$
by $W_{t}(\bolds{\eta})$ and $\bar{W}_{t}(\bolds{\eta})$,
respectively.
By making use of
\begin{eqnarray*}
\inf_{\bolds{\eta}\in B_{\tau_1}(\bolds{\eta}_0)}
\lambda_{\min}\Biggl(n^{-1}\sum_{t=2}^n W_{t}(\bolds{\eta})\Biggr)
&\geq&\inf_{\bolds{\eta}\in B_{\tau_1}(\bolds{\eta }_0)}\lambda_{\min}\Biggl(
n^{-1}\sum_{t=2}^n \bar{W}_{t}(\bolds{\eta})\Biggr)\\
&&{} - \sup_{\bolds{\eta}\in
B_{\tau_1}(\bolds{\eta}_0)} \Biggl\|n^{-1}\sum_{t=2}^{n}\bigl[W_{t}(\bolds{\eta})-
\bar{W}_{t}(\bolds{\eta})\bigr]\Biggr\|,
\end{eqnarray*}
(B.6) of \citet{ChaIng11},
(\ref{222})--(\ref{224}) and (\ref{26}),
we get from an argument similar to that used to prove (\ref{215})
that for $\bar{M}>2/\underline{c}$,
\[
P\Biggl\{\sup_{\bolds{\eta}\in B_{\tau_1}(\bolds{\eta}_0)}
\lambda_{\min}^{-1}\Biggl(n^{-1}\sum_{t=2}^n W_{t}(\bolds{\eta})\Biggr)
>\bar{M}\Biggr\}=o\bigl(n^{-q}\bigr),
\]
which gives (\ref{213}). As the proof of (\ref{214}) is similar to
(\ref{213}), details are omitted.

To prove (\ref{210}), first note that by (\ref{26}) and (B.5) of \citet{ChaIng11}, 
we have for any
$1\le i,j \le\bar{p}$,
%
%
\begin{eqnarray}
\label{225} && \mathrm{E}\Biggl\{\sup_{\bolds{\eta}\in B_{\tau
_1}(\bolds{\eta}_0)} \Biggl|n^{-1/2}
\sum_{t=2}^n\varepsilon_t
\nabla^2\varepsilon_t(\bolds{\eta})_{i,j}\Biggr|^{q_1}
\Biggr\}
\nonumber\\[-8pt]\\[-8pt]
&&\qquad\le C \Biggl[\Biggl\{\sum_{s=1}^{n-1}c^{2}_{s, ij}(
\bolds{\eta}_0) \Biggr\}^{q_{1}/2}+ \Biggl\{\sum
_{s=1}^{n-1}A_{s}(i,j)\Biggr
\}^{q_{1}/2}\Biggr].
\nonumber
\end{eqnarray}
Combining (\ref{225}) with (\ref{220}) gives the desired conclusion.
Finally, by (\ref{26}), (\ref{218}), (\ref{219}), (\ref{222}), (\ref
{223}), Lemma 2 of \citet{Wei87} 
and an argument used in the proof of (\ref{A12}) in Appendix A, we
have for any $1\leq, i, j\leq\bar{p}$,
%
%
\begin{eqnarray}
\label{226} && \mathrm{E}\Bigl\{\sup_{\bolds{\eta}\in B_{\tau
_1}(\bolds{\eta}_0)} \bigl|
\nabla^2\varepsilon_t(\bolds{\eta})_{i,j}\bigr|^{4q_1}
\Bigr\}
\nonumber\\
&&\qquad\leq C \Biggl[ \Biggl\{\sum_{s=1}^{\infty}
\max_{{\mathbf j}\in{\mathbf J}(m, \bar{p}), 1\leq m \leq\bar{p}} \sup
_{\bolds{\eta}\in B_{\tau_1}(\bolds{\eta}_{0})} \bigl({\mathbf
D_j}c_{s,ij}(\bolds{\eta})\bigr)^2 \Biggr
\}^{2q_{1}}\\
&&\qquad\quad\hspace*{126.5pt}{} + \Biggl\{\sum_{s=1}^{\infty}
c^{2}_{s,ij}(\bolds{\eta}_{0})\Biggr
\}^{2q_{1}} \Biggr]
\leq C
\nonumber
\end{eqnarray}
and
%
%
\begin{eqnarray}
\label{227} && \mathrm{E}\Bigl\{\sup_{\bolds{\eta}\in B_{\tau
_1}(\bolds{\eta}_0)} \bigl|\nabla
\varepsilon_t(\bolds{\eta})_{i}\bigr|^{4q_1}
\Bigr\}
\nonumber\\
&&\qquad\leq C \Biggl[ \Biggl\{\sum_{s=1}^{\infty}
\max_{{\mathbf j}\in{\mathbf J}(m, \bar{p}), 1\leq m \leq\bar{p}} \sup
_{\bolds{\eta}\in B_{\tau_1}(\bolds{\eta}_{0})} \bigl(\mathbf
{D}_jb_{s,i}(\bolds{\eta})\bigr)^2 \Biggr
\}^{2q_{1}} \\
&&\qquad\quad\hspace*{125.5pt}{} + \Biggl\{\sum_{s=1}^{\infty}
b^{2}_{s,i}(\bolds{\eta}_{0})\Biggr
\}^{2q_{1}} \Biggr]
\leq C,
\nonumber
\end{eqnarray}
and hence (\ref{211}) and (\ref{212}) hold.
This completes the proof of Theorem~\ref{thm21}.
\end{pf}

We close this section with a brief discussion of generalizing (\ref
{16}) to the linear process
%
%
\begin{equation}
\label{228} y_{t}=m_{t}(\bolds
\eta_{0})+\varepsilon_{t},
\end{equation}
where $\varepsilon_{t}$'s obey (A1),
$\bolds\eta_{0}=(\bolds\theta_0, d_{0})^{\mathsf{T}}$ is an
unknown $\bar{p}$-dimensional vector
with $d_{0} \in D$ and $\bolds\theta_0$ lying in a given compact set
$V \subset R^{\bar{p}-1}$, and
$m_t(\bolds{\eta})=m_{t}(\bolds{\eta}, y_{t-1},\ldots, y_{1})$
admits a linear representation $\sum_{s=1}^{t-1}\tilde{c}_s(\bolds
{\eta})\varepsilon_{t-s}$
with $\tilde{c}_s(\bolds{\eta})$'s being twice differentiable on
$V \times D$.
Assume that $\bolds\eta_{0} \in\operatorname{int} V \times D$
and $\tilde{c}_s(\bolds{\eta})$'s satisfy
some identifiability conditions leading to (\ref{A22}) in the \hyperref
[app]{Appendix}
and (0.1) in the supplementary material [\citet{ChaHuaIng}],
and some smoothness conditions similar to (\ref{218}), (\ref{219}),
(\ref{222}), (\ref{223})
and (\ref{A12}). Then
the same argument used in the proof of Theorem~\ref{thm21} shows that
(\ref{16}) is still valid under (\ref{228}).
Note that these identifiability and smoothness conditions are readily
fulfilled not only by (\ref{11}),
but also by (\ref{11}) with the ARMA component being replaced by the
exponential-spectrum model of
\citet{Blo73}. Moreover, when the ARMA component of (\ref{11}) is
replaced by the more general one given in (1.3)
of \citet{HuaRob11},
these conditions can also be ensured by their assumptions A1 and A3,
with A1(ii), A2(ii) and A2(iii) suitably modified.
\section{Mean squared prediction errors}\label{sec3}
\label{section3}
One important and intriguing application of Theorem~\ref{thm21} is the
analysis of
mean squared prediction errors. Assume that $y_{1},\ldots, y_{n}$ are
generated by model (\ref{11}).
To predict $y_{n+h}$,\break $h\geq1$, based on $y_{1},\ldots, y_{n}$, we
first adopt the one-step CSS predictor,\break
$\hat{y}_{n+1}(\hat{\bolds{\eta}}_{n})=y_{n+1}-\varepsilon
_{n+1}(\hat{\bolds{\eta}}_{n})$, to forecast $y_{n+1}$,
noting that $\hat{y}_{n+1}(\hat{\bolds{\eta}}_{n})$ depends solely
on $y_{1},\ldots, y_{n}$.
Define $p_{t}(\bolds{\eta}, y_{t-1},\ldots,
y_{1}):=y_{t}-\varepsilon_{t}(\bolds{\eta})=(1-(1-\break B)^{d}A_{1,
\bolds{\theta}}(B)A^{-1}_{2, \bolds{\theta}}(B))y_{t}$.
Then $y_{n+h}$, $h \geq2$, can be predicted recursively by the
$h$-step CSS predictor,
%
%
\begin{equation}
\label{31} \hat{y}_{n+h}(\hat{\bolds{\eta}}_{n}):=p_{n+h}
\bigl(\hat{\bolds{\eta}}_{n}, \hat{y}_{n+h-1}(\hat{
\bolds{\eta}}_{n}),\ldots, \hat{y}_{n+1}(\hat{
\bolds{\eta}}_{n}), y_{n},\ldots, y_{1}
\bigr).
\end{equation}
When restricted to the short-memory AR case where $p_{t}(\bolds
{\eta})=(1-A_{1, \bolds{\theta}}(B))y_{t}$,
$\hat{y}_{n+h}(\hat{\bolds{\eta}}_{n})$ is called the plug-in
predictor in\vadjust{\goodbreak} \citet{Ing03}. 
Sections~\ref{sec3.1} and~\ref{sec3.2} provide an asymptotic expression
for the MSPE of
$\hat{y}_{n+h}(\hat{\bolds{\eta}}_{n})$,
$\mathrm{E}\{y_{n+h}-\hat{y}_{n+h}(\hat{\bolds{\eta}}_n)\}^2$,
with $h=1$ and $h>1$, respectively.

\subsection{One-step prediction}\label{sec3.1}
\label{sec31}
In this section, we apply Theorem~\ref{thm21}
to analyze $\mathrm{E}\{y_{n+1}-\hat{y}_{n+1}(\hat{\bolds{\eta
}}_n)\}^2$.
In particular, it is shown in Theorem~\ref{thm31} that
the contribution of the estimated parameters to the one-step MSPE,
$ \mathrm{E}\{y_{n+1}-\hat{y}_{n+1}(\hat{\bolds{\eta}}_n)\}
^2-\sigma^2$,
is proportional to the number of parameters.
We start with the following auxiliary lemma.

%
\begin{lemma}
\label{lemma31}
Assume (\ref{11})--(\ref{15}), \textup{(A1)} and
%
%
\begin{equation}
\label{32} \sup_{t\ge1}\mathrm{E}|\varepsilon_t|^{\gamma}<
\infty\qquad\mbox{for some } \gamma>4.
\end{equation}
Then, for any $\delta>0$ such that
$\Pi\times D-B_{\delta}(\bolds{\eta}_0)\ne\varnothing$,
%
%
\begin{equation}
\label{33} n\mathrm{E}\bigl[\bigl\{\varepsilon_{n+1}(\hat{\bolds{
\eta}}_n)-\varepsilon_{n+1}(\bolds{
\eta}_0)\bigr\}^2 I_{\{\hat{\bolds{\eta}}_n\in\Pi\times D-B_{\delta
}(\bolds{\eta}_0)\}}\bigr]=o(1).
\end{equation}
\end{lemma}

\begin{pf}
Let $4<\gamma_{1}<\gamma$ and $0< v<(\gamma_{1}-4)/(2\gamma_{1}+8)$.
Also let $B_{j,v},\allowbreak j \geq0$, be defined as in Lemma~\ref{lemma2} and
$W$ be
defined as in the proof of Corollary~\ref{cor21}. By
Cauchy--Schwarz's inequality, the left-hand side of (\ref{33}) is
bounded above by
%
%
\begin{equation}
\label{34} n\sum_{j=0}^{W}\mathrm{E}^{1/2} \Bigl\{\sup_{\bolds{\eta}\in
B_{j,v}} \bigl(
\varepsilon_{n+1}(\bolds{\eta})-\varepsilon_{n+1}(
\bolds{\eta}_0)\bigr)^{4} \Bigr\} P^{1/2}(
\hat{\bolds{\eta}}_n \in B_{j,v}).
\end{equation}
By the compactness of $B_{j, v}$, (\ref{A39}) and (\ref{A40}) in
the \hyperref[app]{Appendix}
and an argument similar to that used to prove (\ref{227}), it follows that
%
%
\begin{eqnarray}
\label{35} &&\mathrm{E}^{1/2} \Bigl\{\sup_{\bolds{\eta}\in B_{j,v}}
\bigl(
\varepsilon_{n+1}(\bolds{\eta})-\varepsilon_{n+1}(
\bolds{\eta}_0)\bigr)^{4} \Bigr\}
\nonumber
\\
&&\qquad\leq C \Biggl\{\sup_{\bolds{\eta}\in B_{j,v}}\sum_{s=1}^{n}b^2_s(
\bolds{\eta}) +\sum_{s=1}^{n}\max
_{{\mathbf j}\in J(m,\bar{p}),1\le m\le\bar{p}} \sup_{\bolds{\eta}\in
B_{j,v}}\bigl(D_{\mathbf j}b_s(
\bolds{\eta})\bigr)^2 \Biggr\}
\\
&&\qquad= \cases{ O\bigl((\log n)^3\bigr), &\quad$j=0$,
\vspace*{2pt}\cr
O
\bigl(n^{2vj}(\log n)^2\bigr), &\quad$j\ge1$.}
\nonumber
\end{eqnarray}
Moreover, Lemma~\ref{lemma2} yields that for any $\theta>0$,
%
%
\begin{eqnarray}
\label{36} P^{1/2}(\hat{\bolds{\eta}}_n \in
B_{j,v}) &\leq& C \biggl\{\mathrm{E} \biggl(\frac{{\sup_{\bolds{\eta}\in B_{j,v}}}
|{\sum_{t=1}^n}(\varepsilon_t(\bolds{\eta})-\varepsilon(\bolds
{\eta}_0))\varepsilon_t| } {
\inf_{\bolds{\eta}\in B_{j,v}}\sum_{t=1}^n
(\varepsilon_t(\bolds{\eta})-\varepsilon_t(\bolds{\eta}_0))^2}
\biggr)^{\gamma_1} \biggr\}^{1/2}
\nonumber\\[-8pt]\\[-8pt]
&=& \cases{\displaystyle O \biggl(\biggl(\frac{(\log n)^{3/2}}{n^{1/2}}\biggr
)^{\gamma
_{1}/2}(\log
n)^\theta\biggr), &\quad$j=0$,
\vspace*{2pt}\cr
\displaystyle O \biggl(\biggl(\frac{\log n}{n^{1/2+jv-2v}}
\biggr)^{\gamma_{1}/2}(\log n)^\theta\biggr), &\quad$j\ge1$.}\nonumber
\end{eqnarray}
Combining (\ref{34})--(\ref{36}), we obtain for some $s>0$,
$n\mathrm{E}[\{\varepsilon_{n+1}(\hat{\bolds{\eta}}_n)-\break\varepsilon
_{n+1}(\bolds{\eta}_0)\}^2 \times
I_{\{\hat{\bolds{\eta}}_n\in\Pi\times D-B_{\delta}(\bolds{\eta
}_0)\}}]
=O(n^{1+2v-(\gamma_{1}/2)(1/2-v)}(\log n)^{s})=o(1)$,\break
where the last equality is ensured by the ranges of $\gamma_{1}$ and
$v$ given above.
As a result, (\ref{33}) is proved.
\end{pf}

Equipped with Lemma~\ref{lemma31}, we are now in a position to state
and prove Theorem~\ref{thm31}.
%
%
\begin{theorem}
\label{thm31}
Suppose that the assumptions of Theorem~\ref{thm21} hold except that
(\ref{26}) is replaced by
%
%
\begin{equation}
\label{101} \sup_{t\ge1}\mathrm{E}|\varepsilon_t|^{\gamma}<
\infty\qquad\mbox{for some } \gamma>10.
\end{equation}
Assume also that $\varepsilon_{t}$'s are i.i.d. random variables.
Then
%
%
\begin{equation}
\label{37} \lim_{n\rightarrow\infty} n\bigl[ \mathrm{E}\bigl
\{y_{n+1}-\hat{y}_{n+1}(\hat{\bolds{\eta}}_n)
\bigr\}^2 - \sigma^2 \bigr] =\bar{p}\sigma^2.
\end{equation}
\end{theorem}

\begin{pf}
Let $0<\tau_1<1/2$ satisfy $B_{\tau_{1}}(\bolds{\eta}_{0}) \subset
\Pi\times D$.
Define ${\cal D}_{n}=\{\hat{\bolds{\eta}}_{n} \in B_{\tau
_{1}}(\bolds{\eta}_{0})\}$
and ${\cal D}^{c}_n=\{\hat{\bolds{\eta}}_{n} \in\Pi\times
D-B_{\tau_{1}}(\bolds{\eta}_{0})\}$.
By Taylor's theorem,
%
%
\begin{eqnarray}
\label{38}
&&
n^{1/2}\bigl(y_{n+1}-\hat{y}_{n+1}(
\hat{\bolds{\eta}}_n)-\varepsilon_{n+1}\bigr)\nonumber\\
&&\qquad=
n^{1/2}\bigl(\nabla\varepsilon_t(\bolds{
\eta}_0)\bigr)^{\mathsf{T}}(\hat{\bolds{\eta}}_n-
\bolds{\eta}_0)I_{{\cal D}_n}
\nonumber\\[-8pt]\\[-8pt]
&&\qquad\quad{}+\frac{n^{1/2}}{2}(\hat{\bolds{\eta}}_n-\bolds{
\eta}_0)^{\mathsf{T}} \nabla^2\varepsilon_{n+1}
\bigl(\bolds{\eta}^\ast\bigr) (\hat{\bolds{
\eta}}_n-\bolds{\eta}_0)I_{{\cal D}_n}\nonumber\\
&&\qquad\quad{}+n^{1/2}\bigl(\varepsilon_{n+1}(\hat{\bolds{
\eta}}_n)-\varepsilon_{n+1}(\bolds{
\eta}_0)\bigr)I_{{\cal D}^{c}_n},
\nonumber
\end{eqnarray}
where $\|\bolds{\eta}^\ast-\bolds{\eta}_0\|\le\|\hat
{\bolds{\eta}}_n-\bolds{\eta}_0\|$.
Since Lemma~\ref{lemma31} ensures that the second moment of the third
term on the right-hand side of
(\ref{38}) converges to 0, the desired conclusion (\ref{37}) follows
immediately from
%
%
\begin{equation}
\label{39} \lim_{n\rightarrow\infty}\mathrm{E}\bigl[n\bigl\{\bigl
(\nabla
\varepsilon_{n+1}(\bolds{\eta}_0)
\bigr)^{\mathsf{T}} (\hat{\bolds{\eta}}_n-\bolds{
\eta}_0)\bigr\}^2I_{{\cal D}_n}\bigr] =\bar{p}
\sigma^2
\end{equation}
and
%
%
\begin{equation}
\label{310} \mathrm{E}\bigl[ n\bigl\{(\hat{\bolds{\eta}}_n-
\bolds{\eta}_0)^{\mathsf{T}} \nabla^{2}
\varepsilon_{n+1}\bigl(\bolds{\eta}^{*}\bigr) (\hat{
\bolds{\eta}}_n-\bolds{\eta}_0)\bigr
\}^{2}I_{{\cal D}_{n}}\bigr] =o(1).
\end{equation}

Note first that by Theorem 2.2 of
\citet{HuaRob11}, 
%
%
\begin{equation}
\label{311} n^{1/2}(\hat{\bolds{\eta}}_n-\bolds{
\eta}_0) \Rightarrow\mathbf{Q},
\end{equation}
where
${\mathbf Q}$ is distributed as $N({\mathbf0}$, $\sigma^{2}\Gamma
^{-1}(\bolds{\eta}_{0}))$
with $\Gamma(\bolds{\eta}_{0})_{i,j}=\break
\lim_{t \to\infty}\mathrm{E}(\nabla
\varepsilon_{t}(\bolds{\eta}_{0})_{i}\*
\nabla\varepsilon_{t}(\bolds{\eta}_{0})_{j})=\sigma^{2}$
$\sum_{s=1}^{\infty}b_{s,i}(\bolds{\eta}_{0})b_{s,i}(\bolds
{\eta}_{0})$
for $1\leq i, j\leq\bar{p}$ and $\Rightarrow$ denotes convergence in
distribution.
[Note that
$\Gamma(\bolds{\eta}_{0})$
is independent of $d_{0}$ and
$\Gamma(\bolds{\eta}_{0})_{\bar{p},\bar{p}}=\pi^{2}\sigma^{2}/6$.]
Define
$\nabla\varepsilon_{n+1,m}(\bolds{\eta}_0)
=(\sum_{s=1}^m b_{s, i}(\bolds{\eta}_0)\varepsilon_{n+1-s})_{1\leq
i\leq\bar{p}}$.
Then by (\ref{311}) and the independence between $\nabla\varepsilon
_{n+1,m}(\bolds{\eta}_0)$
and\break $n^{1/2}(\hat{\bolds{\eta}}_{n-m}-\bolds{\eta}_0)$,
%
%
\begin{equation}
\label{312}\quad Z_{n,m}=n^{1/2}\bigl(\nabla\varepsilon_{n+1,m}(
\bolds{\eta}_0)\bigr)^{\mathsf{T}}(\hat{\bolds{
\eta}}_{n-m}-\bolds{\eta}_0) \Rightarrow{\mathbf
F}_m^{\mathsf{T}}{\mathbf Q}\qquad\mbox{as } n\rightarrow\infty
\end{equation}
and
%
%
\begin{equation}
\label{313} {\mathbf F}^{\mathsf{T}}_m{\mathbf Q}\Rightarrow{\mathbf
F}^{\mathsf{T}}{\mathbf Q}\qquad\mbox{as } m\rightarrow\infty,
\end{equation}
where ${\mathbf F}$
and ${\mathbf F}_m$, independent of ${\mathbf Q}$, have the same distribution
as those of
$(\sum_{s=1}^\infty b_{s,i}(\bolds{\eta}_0)\varepsilon_s,)_{1\leq
i \leq\bar{p}}$
and $\nabla\varepsilon_{m+1,m}(\bolds{\eta}_0)$, respectively.
By making use of
(\ref{213}), (\ref{227}), (\ref{311}), $\hat{\bolds{\eta
}}_{n-m} \to_{p} \bolds{\eta}_{0}$ as $n \to\infty$,
and $\nabla S_{n-m}(\hat{\bolds{\eta}}_{n-m})=\mathbf{0}$ on
$\{\hat{\bolds{\eta}}_{n-m} \in B_{\tau_{1}}(\bolds{\eta
}_{0})\}$,
we obtain that for any $\epsilon>0$,
\[
\lim_{m \to\infty} \limsup_{n \to\infty} P\bigl
\{\bigl|n^{1/2}\bigl(\nabla\varepsilon_{n+1}(\bolds{
\eta}_0)\bigr)^{\mathsf{T}}(\hat{\bolds{\eta}}_n-
\bolds{\eta}_0)-Z_{n,m}\bigr|>\epsilon\bigr\}=0,
\]
which, together with $\lim_{n \to\infty}P({\cal D}_{n})=1$, Theorem
4.2 of \citet{Bil68}, 
(\ref{312}), (\ref{313})
and the continuous mapping theorem, yields
%
%
\begin{equation}
\label{314} n\bigl\{\bigl(\nabla\varepsilon_{n+1}(\bolds{
\eta}_0)\bigr)^{\mathsf{T}}(\hat{\bolds{\eta}}_n-
\bolds{\eta}_0)\bigr\}^2I_{{\cal D}_n}
\Rightarrow\bigl({\mathbf F}^{\mathsf{T}}{\mathbf Q}\bigr)^2.
\end{equation}
Let $5<v<\gamma/2$ and $\theta=(\gamma/v)-2$.
It follows from (\ref{101}), Theorem~\ref{thm21} and H\"{o}lder's
inequality that
%
%
\begin{eqnarray}
\label{102}
&&
\mathrm{E}\bigl|\bigl(\nabla\varepsilon_{n+1}(\bolds{
\eta}_0)\bigr)^{\mathsf{T}}n^{1/2} (\hat{\bolds{
\eta}}_n-\bolds{\eta}_0)\bigr|^{2+\theta} \nonumber\\
&&\qquad\le\bigl
\{\mathrm{E}\bigl\|\nabla\varepsilon_{n+1}(\bolds{\eta}_0)
\bigr\|^{\gamma}\bigr\}^{{1}/{v}}
\bigl\{\mathrm{E}\bigl\|n^{1/2}(\hat{\bolds{\theta}}_n-
\bolds{\theta}_0)\bigr\|^{\gamma/(v-1)}\bigr\}^{({v-1})/{v}} \\
&&\qquad=O(1)
\nonumber
\end{eqnarray}
and hence $n\{(\nabla\varepsilon_{n+1}(\bolds{\eta}_0))^{\mathsf
{T}}(\hat{\bolds{\eta}}_n-\bolds{\eta}_0)\}^2I_{{\cal D}_n}$ is
uniformly integrable. This, (\ref{314}) and $\mathrm{E} ({\mathbf
F}^{\mathsf{T}}{\mathbf Q})^{2}= \bar{p} \sigma^{2}$ together imply
(\ref{39}).

On the other hand,
since on ${\cal D}_{n}$, $\|\hat{\bolds{\eta}}_n-\bolds{\eta
}_0\|< \tau_{1}$,
we have for any $0<\theta<2$,
$\|\hat{\bolds{\eta}}_n-\bolds{\eta}_0\|^{4}I_{{\cal D}_{n}}
\leq K \|\hat{\bolds{\eta}}_n-\bolds{\eta}_0\|^{2+\theta
}I_{{\cal D}_{n}}$,
where $K$ is some positive constant depending only on $\theta$ and $\tau_{1}$.
Let $0< \theta< \min\{4^{-1}(\gamma-2)-2, 2\}$.
Then, it follows from Theorem~\ref{thm21}, (\ref{101}) and H\"
{o}lder's inequality that
%
%
\begin{eqnarray}
\label{315}
&&\mathrm{E}\bigl\{ n^{1/2}(\hat{\bolds{
\eta}}_n-\bolds{\eta}_0)^{\mathsf{T}}
\nabla^{2} \varepsilon_{n+1}\bigl(\bolds{
\eta}^{*}\bigr) (\hat{\bolds{\eta}}_n-\bolds{
\eta}_0)I_{{\cal
D}_{n}}\bigr\}^{2}
\nonumber
\\
&&\qquad\le K \mathrm{E}\bigl( n\|\hat{\bolds{\eta}}_n-\bolds{
\eta}_0\|^{2+\theta}\bigl\|\nabla^{2}
\varepsilon_{n+1}\bigl(\bolds{\eta}^{*}\bigr)
\bigr\|^{2}\bigr)
\nonumber\\
&&\qquad\le K n^{-\theta/2}\bigl(\mathrm{E}\bigl\|n^{1/2}(\hat{\bolds{
\eta}}_n-\bolds{\eta}_0) \bigr\|^{{(2+\theta)\gamma}/({\gamma-2})}
\bigr)^{({\gamma-2})/{\gamma}} \\
&&\qquad\quad{}\times\Bigl\{\mathrm{E}\Bigl(\sup
_{\bolds{\eta} \in B_{\tau_{1}(\bolds{\eta}_{0})}}\bigl\|
\nabla^{2} \varepsilon_{n+1}(\bolds{\eta})
\bigr\|^{\gamma}\Bigr)\Bigr\}^{{2}/{\gamma}}
\nonumber
\\
&&\qquad= O\bigl(n^{-\theta/2}\bigr)\Bigl\{\mathrm{E}\Bigl(\sup
_{\bolds
{\eta} \in B_{\tau_{1}(\bolds{\eta}_{0})}
}
\bigl\|\nabla^{2} \varepsilon_{n+1}(\bolds{\eta})
\bigr\|^{\gamma}\Bigr)\Bigr\}^{{2}/{\gamma}}.
\nonumber
\end{eqnarray}
An argument similar to that used to prove (\ref{226}) also yields
that the expectation on the right-hand side of (\ref{315}) is bounded
above by a finite constant, and hence
(\ref{310}) holds true. This completes the proof of the theorem.
\end{pf}
Theorem~\ref{thm31} asserts that the second-order MSPE of
$\hat{y}_{n+1}(\hat{\bolds{\eta}}_n)$,
$\bar{p}\sigma^{2}n^{-1}+o(n^{-1})$, only depending the number of
estimated parameters, has nothing to do with dependent structure of
the underlying process. This result is particularly interesting when
compared with the second-order MSPE of the LS predictor in
integrated AR models. To see this, assume first that
there is a forecaster who believes that the true model is possibly an
integrated AR($p_{1}$) model,
%
%
\begin{eqnarray}
\label{316}
&&\bigl(1-\tilde{\alpha}_{1}B-\cdots-\tilde{
\alpha}_{p_{1}}B^{p_{1}}\bigr)y_{t}\nonumber\\[-8pt]\\[-8pt]
&&\qquad=(1-B)^{v_{0}}\bigl(1-\theta_{1}B-\cdots\theta_{p_{1}-v_{0}}B^{p_{1}-v_{0}}
\bigr)y_{t}=\varepsilon_{t},\nonumber
\end{eqnarray}
where $v_{0} \in\{0, 1,\ldots, p_{1}\}$ is unknown and $1-\theta
_{1}z-\cdots\theta_{p_{1}-d}z^{p_{1}-v_{0}} \neq0$ for all $|z|\leq
1$. Then
it is natural for this forecaster to predict $y_{n+1}$ using
the LS predictor $\tilde{y}_{n+1}$, in which
$\tilde{y}_{n+1}=\mathbf{y}^{\mathsf{T}}_{n}(p_{1})\tilde{\bolds
{\alpha}}_{n}(p_{1})$
with $\mathbf{y}^{\mathsf{T}}_{t}(p_{1})=(y_{t},\ldots, y_{t-p_{1}+1})$
and $\tilde{\bolds{\alpha}}_{n}(p_{1})$ satisfies $\sum
_{t=p_{1}}^{n-1}\mathbf{y}_{t}(p_{1})\mathbf{y}^{\mathsf
{T}}_{t}(p_{1})\tilde{\bolds{\alpha}}_{n}(p_{1})=
\sum_{t=p_{1}}^{n-1}\mathbf{y}_{t}(p_{1})y_{t+1}$.
On the other hand, another forecaster who
has doubts on whether the $v_{0}$ in (\ref{316}) is really an integer,
chooses a more flexible alternative as follows:
%
%
\begin{equation}
\label{317} \bigl(1-\alpha_{1}B-\cdots-\alpha_{p_{1}}B^{p_{1}}
\bigr) (1-B)^{d_{0}}y_{t}=\varepsilon_{t},
\end{equation}
where $L_{1}\leq0 \leq d_{0} \leq p_{1} \leq U_{1}$ with $-\infty
<L_{1}<U_{1}<\infty$ being some prescribed numbers,
and $1-\sum_{j=1}^{p_{1}}\alpha_{j}z^{j}$ satisfies (\ref{13}).
Clearly, model (\ref{317}), including model (\ref{316}) as a
particular case, is itself a special case of model (\ref{11}) with $p_{2}=0$,
and hence the CSS predictor, $\hat{y}_{n+1}(\hat{\bolds{\eta
}}_{n})$, obtained from (\ref{317})
is adopted naturally by the second forecaster.

If the data are truly generated by (\ref{316}), then Theorem 2 of \citet{IngSinYu10} 
shows that
under certain regularity conditions,
%
%
\begin{equation}
\label{319} \lim_{n\rightarrow\infty} n\bigl[ \mathrm{E}
\{y_{n+1}-\tilde{y}_{n+1}\}^2 -
\sigma^2 \bigr] =\bigl(p_{1}+v^{2}_{0}
\bigr)\sigma^2.
\end{equation}
In addition, by Theorem~\ref{thm31} (which is still valid in the case
of $p_{2}=0$), we have
%
%
\begin{equation}
\label{320} \lim_{n\rightarrow\infty} n\bigl[ \mathrm{E}\bigl
\{y_{n+1}-\hat{y}_{n+1}(\hat{\bolds{\eta}}_n)
\bigr\}^2 - \sigma^2 \bigr] =(p_{1}+1)
\sigma^2.
\end{equation}
As shown in (\ref{319}) and (\ref{320}), while the second-order MSPE
of the LS predictor $\tilde{y}_{n+1}$ increases as the strength of dependence
in the data does (i.e., $v_{0}$ increases), the second-order MSPE of
the CSS predictor $\hat{y}_{n+1}(\hat{\bolds{\eta}}_n)$ does not
vary with $v_{0}$.
These equalities further indicate the somewhat surprising fact that
for an integrated AR model, even the most popular LS predictor can
be inferior to the CSS predictor, if the integration order is large.
To further illustrate (\ref{319}) and (\ref{320}), we conduct a
simulation study to compare the empirical estimates of $n[
\mathrm{E}\{y_{n+1}-\tilde{y}_{n+1}\}^2 - \sigma^2]$ and $n[
\mathrm{E}\{y_{n+1}-\hat{y}_{n+1}(\hat{\bolds{\eta}}_n)\}^2 -
\sigma^2 ]$ for (\ref{316}) with $p_{1}=3$ and $v_{0}=0,1$ and 2.
These estimates, obtained based on 5000 replications for $n =
1000$, are summarized in Table~\ref{table1}. As observed in Table
\ref{table1}, the empirical estimates of
%
%
\begin{table}
\tablewidth=275pt
\caption{The empirical estimates of the second-order MSPEs of the CSS
predictor (with $p_{1}=3$ and $p_{2}=0$)
and the LS predictors with $p_{1}=3$}
\label{table1}
\begin{tabular*}{\tablewidth}{@{\extracolsep{\fill}}lcc@{}}
\hline
\textbf{True model} & \textbf{CSS predictor} & \textbf{LS predictor} \\
\hline
$(1+0.5B)(1-B)^2y_t=\varepsilon_t$
&4.0689 &6.8409 \\
$(1-0.25B^2)(1-B)y_t=\varepsilon_t$
&4.3732 &4.1975 \\
$(1-0.2B-0.25B^2+0.5B^3)y_t=\varepsilon_t$
&4.1828 &3.1686 \\
\hline
\end{tabular*}
\end{table}
$n[\mathrm{E}\{y_{n+1}-\hat{y}_{n+1}(\hat{\bolds{\eta}}_n)\}^2 -
\sigma^2 ]$ are quite close to 4 for all three models, whereas those
of $n[ \mathrm{E}\{y_{n+1}-\tilde{y}_{n+1}\}^2 - \sigma^2]$ are not
distant from 7, 4 and 3 for $v_{0}=2, 1$ and 0, respectively. Hence
all these estimates align with their corresponding limiting values
given in (\ref{319}) and (\ref{320}). This ``dependency-free''
feature of the CSS predictor in the one-step case, however, vanishes
in the multi-step case, as will be seen in the next section.


\subsection{Multi-step prediction}\label{sec3.2}
\label{sec33}
Note that under (\ref{11}),
$y_{t}=\sum_{j=0}^{t-1}\underline{c}_{s}(\bolds{\eta
}_{0})\varepsilon_{t-s}$,
where for $\bolds{\eta}= (\bolds{\theta}^{\mathsf{T}},
d)^{\mathsf{T}}
=(\alpha_{1},\ldots, \alpha_{p_{1}},\beta_{1},\ldots, \beta_{p_{2}},
d)^{\mathsf{T}} \in\Pi\times D$,
$\underline{c}_{0}(\bolds{\eta})=1$
and $\underline{c}_{s}(\bolds{\eta})$'s, $s\geq0$, satisfy
$\sum_{s=0}^{\infty}\underline{c}_{s}(\bolds{\eta})z^{s}=A_{2,
\bolds{\theta}}(z)A^{-1}_{1, \bolds{\theta}}(z)(1-z)^{-d}, |z|<1$.
In addition, let $\bar{c}_{0}(d)=1$ and
$\bar{c}_{s}(d)$'s, $s\geq0$, satisfy
$\sum_{s=0}^{\infty}\bar{c}_{s}(d)z^{s}=(1-z)^{-d},\break |z|<1$.
With $v_{t}(d)=(1-B)^{d}y_{t}$, define
$\mathbf{u}_{n}(\bolds{\eta})=(-v_{n}(d),\ldots,
-v_{n-p_{1}+1}(d),\break\varepsilon_{n}(\bolds{\eta}),\ldots,\varepsilon
_{n-p_{2}+1}(\bolds{\eta}))^{\mathsf{T}}$.
Now the $h$-step CSS predictor of $y_{n+h}$ is given by
$\hat{y}_{n+h}(\hat{\bolds{\eta}}_{n})=G_{\hat{d}_{n}}(B)y_{n}+\sum
_{s=0}^{h-1}\bar{c}_{s}(\hat{d}_{n})\hat{v}_{n+h-s}$,
where $G_{d}(B)=B^{-h}-\break(1- B)^{d}\sum_{k=0}^{h-1}\bar
{c}_{k}(d)\*B^{k-h}=(1-B)^{d}\sum_{k=h}^{\infty}\bar{c}_{k}(d)B^{k-h}$,
$\hat{v}_{n+l}=-\mathbf{u}^{\mathsf{T}}_{n} (\hat{\bolds{\eta
}}_{n}) \times\break
A^{l-1}(\hat{\bolds{\theta}}_{n}) \hat{\bolds{\theta}}_{n}$
and
\[
A(\bolds{\theta})= \lleft(
\begin{array} {c|c|c}
\bolds\alpha& \displaystyle \frac{I_{p_1-1}}{{\mathbf0}^{\mathsf{T}}_{p_1-1}} &
{\mathbf0}_{p_1\times p_2}
\\
\hline
\bolds\beta&{\mathbf0}_{p_2\times(p_1-1)} & {\mathbf0}_{p_2}\bigg|
\displaystyle \frac{I_{p_2-1}}{{\mathbf0}^{\mathsf{T}}_{p_2-1}}
\end{array}
\rright).
\]
Here $\bolds{\alpha}=(\alpha_{1},\ldots, \alpha_{p_{1}})^{\mathsf{T}}$,
$\bolds{\beta}=(\beta_{1},\ldots, \beta_{p_{2}})^{\mathsf{T}}$, and
${\mathbf0}_m$, ${\mathbf0}_{m\times n}$ and $I_m$, respectively,
denote the
$m$-dimensional zero vector,
the $m\times n$ zero matrix and the $m$-dimensional identity matrix.
Define $\tilde{v}_{n+l}(\bolds{\eta}_{0})=-\mathbf{u}^{\mathsf{T}}_{n}
(\bolds{\eta}_{0})A^{l-1}(\bolds{\theta}_{0})\bolds
{\theta}_{0}$.
Then, it follows that
$y_{n+h}=G_{d_{0}}(B)y_{n}+\sum_{s=0}^{h-1}\bar
{c}_{s}(d_{0})v_{n+h-s}(d_{0})=G_{d_{0}}(B)y_{n}+ \sum_{s=0}^{h-1}\bar
{c}_{s}(d_{0})
\tilde{v}_{n+h-s}(\bolds{\eta}_{0})+\sum_{s=0}^{h-1}\underline
{c}_{s}(\bolds{\eta}_{0})\varepsilon_{n+h-s}$.
In this section, we establish
an asymptotic expression for
%
%
\begin{eqnarray}
\label{113}
&&\mathrm{E}\bigl[y_{n+h}-\hat{y}_{n+h}(\hat{
\bolds{\eta}}_{n})\bigr]^{2} \nonumber\\
&&\qquad= \sigma^{2}_{h}(
\bolds{\eta}_{0}) + \mathrm{E}\Biggl\{G_{\hat{d}_{n}}(B)y_{n}+
\sum_{s=0}^{h-1}\bar{c}_{s}(
\hat{d}_{n})\hat{v}_{n+h-s}
\\
&&\hspace*{55pt}\qquad\quad{}-G_{d_{0}}(B)y_{n}-\sum_{s=0}^{h-1}
\bar{c}_{s}(d_{0})\tilde{v}_{n+h-s}(\bolds{
\eta}_{0})\Biggr\}^{2},
\nonumber
\end{eqnarray}
where $\sigma^{2}_{h}(\bolds{\eta}_{0})=\sigma^{2}\sum
_{s=0}^{h-1}\underline{c}^{2}_{s}(\bolds{\eta}_{0})$.
To state our result, first express $\Gamma(\bolds{\eta}_{0})$ as
\[
\Gamma(\bolds{\eta}_{0})=\pmatrix{ \Gamma_{11}(
\bolds{\theta}_{0}) &\bolds{\gamma}_{12}(
\bolds{\theta}_{0})
\vspace*{2pt}\cr
\bolds{\gamma}^{\mathsf{T}}_{12}(
\bolds{\theta}_{0}) & \pi^{2}\sigma^{2}/6},
\]
where $\Gamma_{11}(\bolds{\theta}_{0})=(\Gamma(\bolds{\eta
}_{0})_{i,j})_{1\leq i, j\leq p_{1}+p_{2}}$
and $\bolds{\gamma}^{\mathsf{T}}_{12}(\bolds{\theta
}_{0})=(\Gamma(\bolds{\eta}_{0})_{\bar{p},i})_{1\leq i \leq p_{1}+p_{2}}$,
noting that $ \Gamma(\bolds{\eta}_{0})$ is independent of $d_{0}$.
Then
\[
\Gamma^{-1}(\bolds{\eta}_{0})= \pmatrix{ \tilde{
\Gamma}_{11}(\bolds{\theta}_{0}) &\tilde{\bolds{
\gamma}}_{12}(\bolds{\theta}_{0})
\vspace*{2pt}\cr
\tilde{
\bolds{\gamma}}^{\mathsf{T}}_{12}(\bolds{
\theta}_{0}) &\tilde{\gamma}_{22}(\bolds{
\theta}_{0})},
\]
where
$\tilde{\Gamma}_{11}(\bolds{\theta}_{0})=(\Gamma_{11}(\bolds
{\theta}_{0})-\bolds{\gamma}_{12}(\bolds{\theta}_{0})
\bolds{\gamma}^{\mathsf{T}}_{12}(\bolds{\theta}_{0})\gamma
^{-1}_{22}(\bolds{\theta}_{0}))^{-1}$,
$\tilde{\gamma}_{22}(\bolds{\theta}_{0})=(\pi^{2}\sigma^{2}/6-
\bolds{\gamma}^{\mathsf{T}}_{12}(\bolds{\theta}_{0})
\Gamma^{-1}_{11}(\bolds{\theta}_{0}) \bolds{\gamma
}_{12}(\bolds{\theta}_{0}))^{-1}$ and
$\tilde{\bolds{\gamma}}_{12}(\bolds{\theta}_{0})=-
\tilde{\gamma}_{22}(\bolds{\theta}_{0})
\Gamma^{-1}_{11}(\bolds{\theta}_{0}) \bolds{\gamma
}_{12}(\bolds{\theta}_{0})$. Define
$\nabla^{(1)}\varepsilon_{t}(\bolds{\theta}_{0})=
(\nabla\varepsilon_{t}(\bolds{\eta}_{0})_{1},\ldots,
\nabla\varepsilon_{t}(\bolds{\eta}_{0})_{p_{1}+p_{2}})^{\mathsf
{T}}$ [noting that $\nabla\varepsilon_{t}(\bolds{\eta}_{0})_{i},
1\leq i \leq\bar{p}$, is independent of $d_{0}$], $\mathbf{w}_{t,
h}=(\sum_{k=1}^{t-1} \varepsilon_{t-k}/(k+h-1),\ldots,\sum_{k=1}^{t-1}
\varepsilon_{t-k}/k)^\mathsf{T}$,
$Q_{h}(\bolds{\theta}_{0})=\lim_{t \to\infty}\mathrm{E}(\nabla
^{(1)}\varepsilon_{t}(\bolds{\theta}_{0})\mathbf{w}^{\mathsf
{T}}_{t, h})$ and $R(h)=(\gamma_{i,j})_{h\times h}$, in which
$\gamma_{i,i}=6\pi^{-2}\sum_{l=h-i+1}^{\infty}l^{-2}$, $1\leq i \leq
h$,
and $\gamma_{i,j}=\gamma_{j,i}=6\pi^{-2}(j-i)^{-1}\sum
_{l=1}^{j-i}(h-j+l)^{-1}$, $1\le i<j\le h$. Now, an asymptotic
expression for (\ref{113}) is given as follows.

%
\begin{theorem}
\label{thm33}
Under the hypothesis of Theorem~\ref{thm31},
%
%
\begin{eqnarray}
\label{318}
&&
\lim_{n\rightarrow\infty}n \bigl\{ \mathrm{E}
\bigl[y_{n+h}-\hat{y}_{n+h}(\hat{\bolds
\eta}_n)\bigr]^2-\sigma^2_h(
\bolds{\eta}_{0}) \bigr\}\nonumber\\[-8pt]\\[-8pt]
&&\qquad= \bigl\{\underline{f}_{h}(p_{1},
p_{2})+\underline{g}_{h}(\bolds{\eta}_{0})
+ 2J_h(\bolds{\eta}_{0})\bigr\}
\sigma^2,\nonumber
\end{eqnarray}
where\vspace*{1pt}
$\underline{f}_{h}(p_{1}, p_{2})=\operatorname{\mathsf{tr}}\{\Gamma
_{11}(\bolds
{\theta}_{0}) \underline{L}_h(\bolds{\eta}_{0})
\tilde{\Gamma}_{11}(\bolds{\theta}_{0})\underline{L}_h^{\mathsf
{T}}(\bolds{\eta}_{0})\}$,
$\underline{g}_{h}(\bolds{\eta}_{0})=(\pi\sigma^{2}/6)\times\break \tilde{\gamma
}_{22}(\bolds{\theta}_{0})
\underline{\mathbf{c}}^{\mathsf{T}}_{h}(\bolds{\eta
}_{0})R(h)\underline{\mathbf{c}}_{h}(\bolds{\eta}_{0})$, and
$J_h(\bolds{\eta}_{0})=\tilde{\bolds{\gamma}}^{\mathsf
{T}}_{12}(\bolds{\theta}_{0})
\underline{L}_h^{\mathsf{T}}(\bolds{\eta}_{0})
Q_{h}(\bolds{\theta}_{0})\underline{\mathbf{c}}_{h}(\bolds
{\eta}_{0})$,
with
$\underline{\mathbf{c}}_{h}(\bolds{\eta}_{0})=(\underline
{c}_0(\bolds{\eta}_{0}),\ldots,
\underline{c}_{h-1}(\bolds{\eta}_{0}))^{\mathsf{T}}$,
$\underline{L}_h(\bolds{\eta}_{0})=\sum_{s=0}^{h-1}\underline
{c}_s(\bolds{\eta}_{0})$
$\tilde{A}^{h-1-s}(\bolds{\theta}_{0})$ and
\[
\label{115} \tilde{A}(\bolds{\theta})= \lleft(
\begin{array} {c|c} {
\bolds\alpha}\bigg|\displaystyle
\frac{I_{p_1-1}}{{\mathbf0}^{\mathsf{T}}_{p_1-1}} &
{\mathbf0}_{p_1\times p_2}
\\
\hline{\mathbf0}_{p_2\times p_1} & {\bolds\beta}\bigg|\displaystyle\frac
{I_{p_2-1}}{{\mathbf0}^{\mathsf{T}}_{p_2-1}}
\end{array}
\rright).
\]
\end{theorem}
A few comments on Theorem~\ref{thm33} are in order. When $h=1$,
straightforward calculations imply
$\underline{f}_{h}(p_{1},p_{2})+\underline{g}_{h}(\bolds{\eta
}_{0}) +2J_h(\bolds{\eta}_{0})=\bar{p}$,
which leads\vspace*{1pt} immediately to Theorem~\ref{thm31}.
When $p_{1}=p_{2}=0$, $\underline{f}_{h}(p_{1}, p_{2})$
and $2J_h(\bolds{\eta}_{0})$ vanish, and $\tilde{\gamma
}_{22}(\bolds{\eta}_{0})$
and $\underline{\mathbf{c}}_{h}(\bolds{\eta}_{0})$ in
$\underline{g}_{h}(\bolds{\eta}_{0})$ become
$6/(\pi^{2}\sigma^{2})$ and
$\bar{\mathbf{c}}_{h}(d_{0})= (\bar{c}_{0}(d_{0}),\ldots,
\bar{c}_{h-1}(d_{0}))^{\mathsf{T}}$, respectively. As a result, the
right-hand side of (\ref{318}) is simplified to
%
%
\begin{equation}
\label{321} \bar{g}_{h}(d_{0})\sigma^{2}=
\bar{\mathbf{c}}^{\mathsf{T}}_{h}(d_{0})R(h)\bar{
\mathbf{c}}_{h}(d_{0})\sigma^{2},
\end{equation}
yielding the second-order MSPE of the $h$-step CSS predictor for a pure
I($d$) process.
Alternatively, if $d_{0}=0$ is known, then $\underline
{g}_{h}(\bolds{\eta}_{0})$
and $2J_h(\bolds{\eta}_{0})$ vanish, and
the right-hand side of (\ref{318}) becomes
%
%
\begin{equation}
\label{322} \bar{f}_{h}(p_{1}, p_{2})=
\operatorname{\mathsf{tr}}\bigl\{\Gamma_{11}(\bolds{\theta}_{0})
\tilde{L}_h(\bolds{\theta}_{0})
\Gamma^{-1}_{11}(\bolds{\theta}_{0})
\tilde{L}_h^{\mathsf{T}}(\bolds{\theta}_{0})\bigr
\},
\end{equation}
where $\tilde{L}_h(\bolds{\theta}_{0})$ is $ \underline
{L}_h(\bolds{\eta}_{0})$ with $d_{0}=0$.
Note that (\ref{322}) has been obtained by \citet{Yam81} %
under the stationary ARMA($p_{1}, p_{2}$) model through a somewhat
heuristic argument that does not involve the moment bounds of the
estimated parameters.
When $p_{2}=0$, the right-hand side of (\ref{322}) further reduces to
$f_{1, h}(p_{1})$ in (10) of \citet{Ing03}, 
which is the second-order MSPE of the $h$-step plug-in predictor of a
stationary AR($p_{1}$) model. In view\vspace*{-1pt} of the similarity between
$\underline{f}_{h}(p_{1}, p_{2})$ and $\bar{f}_{h}(p_{1}, p_{2})$ and
that between
$\underline{g}_{h}(\bolds{\eta}_{0})$ and $\bar{g}_{h}(d_{0})$,
(\ref{318}) displays not only an interesting structure of the
multi-step prediction formula from the ARMA case to the I($d$) case,
and eventually to the ARFIMA case, but also reveals that the multi-step
MSPE of an ARFIMA model is the sum of one ARMA term,
$\underline{f}_{h}(p_{1}, p_{2})$, one I($d$) term, $\underline
{g}_{h}(\bolds{\eta}_{0})$
and the term $2J_h(\bolds{\eta}_{0})$ that is related to the ARMA
and I($d$) joint effect.
This expression is different from the ones
obtained for the LS predictors in the integrated AR models, in which
the AR and I($d$) joint effect vanishes asymptotically; see Theorem 2.2
of \citet{IngLinYu09} 
and
Theorem 2 of \citet{IngSinYu10} 
for details.

Before leaving this section, we remark that the dependence structure of
(\ref{11}) has a substantial impact on the multi-step MSPE.
To see this, consider the pure I($d$) case.
By (\ref{321}) and a straightforward calculation, it follows that for
any $d_{0} \in R$, there exist
$0<C_{1, d_{0}} \leq C_{2, d_{0}}< \infty$ such that
%
%
\begin{equation}
\label{323} C_{1, d_{0}}h^{-1+2d_{0}} \leq\bar{g}_{h}(d_{0})
\leq C_{2, d_{0}}h^{-1+2d_{0}},
\end{equation}
which shows that for $h>1$, a larger $d_{0}$ (or a stronger dependence
in the
data) tends to result in a larger second-order MSPE.
Finally, if the true model is a random walk model, $y_{t}=\alpha
y_{t-1}+\varepsilon_{t}$, with $\alpha=1$, and is modeled by an
I($d$) process, $(1-B)^{d}y_{t}=\varepsilon_{t}$, in which $d=1$
corresponds to the true model, then by Theorem~\ref{thm33} and (\ref{321}),
$ \lim_{n\rightarrow\infty}n
\{ \mathrm{E}[y_{n+h}-\hat{y}_{n+h}(\hat{d}_n)]^2-h\sigma^{2} \}
=\sigma
^{2}$ for $h=1$, and
the limit is smaller than
$(4.87 h+(1+\log h)^{2}-2(1+\log2h))\sigma^{2}$ for $h \geq2$.
On the other hand, for the $h$-step LS predictor $\tilde{y}_{n+h}$
of the above AR(1) model, it follows from Theorem 2.2 of \citet{IngLinYu09} 
that
$ \lim_{n\rightarrow\infty}n
\{ \mathrm{E}(y_{n+h}-\tilde{y}_{n+h})^2-h\sigma^{2} \}
=2h^{2}\sigma^{2}$,
which is larger than $\sigma^{2}$ when $h=1$,
and larger than $(4.87 h+(1+\log h)^{2}-2(1+\log2h))\sigma^{2}$ when
$h\geq2$.
Hence $\hat{y}_{n+h}(\hat{d}_n)$ is always better than $\tilde
{y}_{n+h}$ in terms of the MSPE.
The convergence rates of their corresponding estimates, however, are
completely reversed because
the LS estimate converges much faster to 1 than
$\hat{d}_{n}$ for a random walk model. This finding is reminiscent of
the fact that when the true model
simultaneously belongs to several different parametric families, the
so-called optimal choice of parametric families
may vary according to different objectives. For a random walk model,
when estimation is the ultimate goal, then LS estimate
may be preferable. On the other hand, for prediction purposes, CSS
predictor is more desirable according to Theorem~\ref{thm33}.

\begin{appendix}\label{app}
\section*{Appendix}

\begin{pf*}{Proof of Lemma~\ref{lemma1}}
We only prove (\ref{21}) for $q\ge1$ because for $0<q<1$, (\ref{21})
is an immediate consequence of for the case of $q\ge1$ and Jensen's inequality.
Since $\Pi\times D-B_{\delta}(\bolds{\eta}_0)$ is compact, there
exists a set of $m$ points
$\{\bolds{\eta}_1,\ldots,\bolds{\eta}_m \}\subset\Pi\times
D-B_{\delta}(\bolds{\eta}_0)$ and a small positive number $0<\delta_1<1$,
depending possibly on $\delta$ and $\Pi$, such that
%
%
\setcounter{equation}{0}
\begin{eqnarray}
\label{A1} \Pi\times D-B_{\delta}(\bolds{\eta}_0)
&\subset&\bigcup_{i=1}^m B_{\delta_1}(
\bolds{\eta}_i),
\\
\label{A2} \|\bolds{\eta}-\bolds{\eta}_0\|&\ge&
\delta/2\quad\mbox{and}\quad\bolds{\theta} \mbox{ obeys
(\ref{13})--(\ref{15})}
\end{eqnarray}
for each $\bolds{\eta}=(\bolds{\theta}^{\mathsf
{T}},d)^{\mathsf{T}}\in\bar{B}_{\delta_1}(\bolds{\eta}_i)$ and
$1\le i\le m$,
where $\bar{B}_{\delta_1}(\bolds{\eta}_i)$ denotes the closure of
$B_{\delta_1}(\bolds{\eta}_i)$.
In view of (\ref{A1}), it suffices for (\ref{21}) to show that
%
%
\begin{eqnarray}
\label{A3}
&&\mathrm{E} \Biggl[ \Biggl\{\inf_{\bolds{\eta}\in\bar
{B}_{\delta_1}(\bolds{\eta}_i)}a_n^{-1}(d)
\sum_{t=1}^n\bigl(\varepsilon_t(
\bolds{\eta})-\varepsilon_t(\bolds{\eta}_0)
\bigr)^2 \Biggr\}^{-q} \Biggr]=O \bigl((\log
n)^{\theta} \bigr),\nonumber\\[-8pt]\\[-8pt]
&&\eqntext{i=1,\ldots,m,}
\end{eqnarray}
hold for any $q\ge1$ and $\theta>0$.
Let $D_i=\{d\dvtx(\bolds{\theta}^{\mathsf{T}},d)^{\mathsf{T}}\in\bar
{B}_{\delta_1}(\bolds{\eta}_i)\}$,
$G_1=\{i\dvtx1\le i\le m, \bar{D}_i=\sup D_i\le d_0-1/2\}$,
$G_2=\{i\dvtx1\le i\le m, \underline{D}_i=\inf D_i\ge d_0-1/2\}$
and $G_3=\{i\dvtx1\le i\le m, \underline{D}_i<d_0-1/2<\bar{D}_i\}$.
Then $\{1,\ldots,m\}=\bigcup_{\ell=1}^3 G_\ell$. In the following, we
first prove (\ref{A3}) for the most challenging case, $i\in G_3$.
The proofs of (\ref{A3}) for the cases of $i\in G_1$ or $G_2$ are
similar, but simpler and are thus omitted.

By the convexity of $x^{-q}, x\ge0$, it follows that for any fixed
$0<\iota<1$,
%
%
\begin{eqnarray}
\label{A4}
&&
\Biggl\{a_n^{-1}(d)\sum
_{i=1}^n\bigl(\varepsilon_t(
\bolds{\eta})-\varepsilon_t(\bolds{\eta}_0)
\bigr)^2 \Biggr\}^{-q} \nonumber\\
&&\qquad\le\Biggl\{n^{-1}\sum
_{t=n\iota+1}^ng_t^2(
\bolds{\eta}) \Biggr\}^{-q}
\\
&&\qquad\le\bigl\{\ell q/(1-\iota)\bigr\}^q z_n^{-1}
\sum_{j=0}^{z_n-1} \Biggl\{\sum
_{r=0}^{\ell q-1}g^2_{n\iota+1+rz_n+j}(
\bolds{\eta}) \Biggr\}^{-q},
\nonumber
\end{eqnarray}
where $\ell>\max\{(2/\alpha_0)+(\ell_1+1)\bar{p}/(\alpha_0q), (\iota
^{-1}-1)/q\}$,
with $\ell_1>2q$, $g_{t}(\bolds{\eta})$ is defined after Lemma \ref
{lemma1},
and $z_n=(1-\iota)n/(\ell q)$.
Here $n\iota$, $\ell q$ and $z_n$ are assumed to be positive integers.
According to (\ref{A4}), if for any $q\ge1$
and all large $n$,
%
%
\begin{eqnarray}
\label{A5}
&&
\mathrm{E} \Biggl[ \Biggl\{\inf_{\bolds{\eta}\in\bar
{B}_{\delta_1}(\bolds{\eta}_i)} \sum
_{r=0}^{\ell q-1}g^2_{n\iota+1+z_nr+j}(
\bolds{\eta}) \Biggr\}^{-q} \Biggr]\le C (\log n)^{5/2},\nonumber\\[-8pt]\\[-8pt]
&&\eqntext{j=0,\ldots,z_n-1,}
\end{eqnarray}
holds, then (\ref{A3}) follows with $\theta=5/2$.
Moreover, since $q$ is arbitrary, this result is easily extended to any
$\theta>0$ using Jensen's inequality.
Consequently, (\ref{A3}) is proved.
In the rest of the proof, we only show that (\ref{A5}) holds for $j=0$
because the proof of (\ref{A5}) for $1\le j\le z_n-1$ is almost identical.
For $j=0$, the left-hand side of (\ref{A5}) is bounded above by
%
%
\begin{eqnarray}
\label{A6}
&&
K+\int_K^\infty P \Biggl(\inf
_{\bolds{\eta}\in\bar{B}_{\delta_1}(\bolds{\eta}_i)}\sum
_{r=0}^{\ell q-1}g^2_{n\iota+1+rz_n}(
\bolds{\eta}) < \mu^{-q^{-1}}, R(\mu) \Biggr) \,d\mu
\nonumber\\
&&\quad{}+ \int_K^\infty P\bigl(R^c(\mu)
\bigr) \,d\mu\\
&&\qquad:= K+\mathrm{(I)}+\mathrm{(II)},
\nonumber
\end{eqnarray}
where $K$, independent of $n$ and not smaller than 1, will be specified
later, with
$c_{\mu}=2\bar{p}^{1/2}\mu^{-(\ell_1+1)/(2q)}$,
\[
R(\mu)=\bigcap_{r=0}^{\ell q-1} \Bigl\{ \mathop{\sup
_{\|\bolds{\eta}_a-\bolds{\eta}_b\|\le c_{\mu}}}_{
\bolds{\eta}_a,\bolds{\eta}_b\in\bar{B}_{\delta_1}(\bolds
{\eta}_i)} \bigl|g_{n\iota+1+rz_n}(\bolds{\eta}_a)-g_{n\iota+1+rz_n}(
\bolds{\eta}_b)\bigr|<2\mu^{-{1}/({2q})} \Bigr\},
\]
and $R^c(\mu)$ is the complement of $R(\mu)$.

We first show that
%
%
\begin{equation}
\label{A9} \mathrm{(II)}\le C(\log n)^{5/2}.
\end{equation}
Define $Q_1^{(i)}=\{ (\bolds{\theta}^{\mathsf{T}},d)^{\mathsf
{T}}\dvtx(\bolds{\theta}^{\mathsf{T}},d)^{\mathsf{T}}\in\bar{B}_{\delta
_1}(\bolds{\eta}_i), d\le d_0-1/2 \}$
and $Q_2^{(i)}=\break\{(\bolds{\theta}^{\mathsf{T}},d)^{\mathsf
{T}}\dvtx(\bolds{\theta}^{\mathsf{T}},d)^{\mathsf{T}}\in\bar{B}_{\delta
_1}(\bolds{\eta}_i), d\ge d_0-1/2\}$.
It is clear that $g_t(\bolds{\eta})$ is differentiable on
$Q_{1,0}^{(i)}$ and $Q_{2,0}^{(i)}$, the interior of
$Q_1^{(i)}$ and $Q_2^{(i)}$, and is continuous on $\bar{B}_{\delta
_1}(\bolds{\eta}_i)$.
By the mean value theorem, we have for any $\bolds{\eta
}_a,\bolds{\eta}_b\in\bar{B}_{\delta_1}(\bolds{\eta}_i)$,
%
%
\begin{eqnarray}
\label{A7}
&& \bigl|g_t(\bolds{\eta}_a)-g_t(
\bolds{\eta}_b)\bigr|\nonumber\\[-8pt]\\[-8pt]
&&\qquad\le\|\bolds{\eta}_a-
\bolds{\eta}_b\| \Bigl( \sup_{\bolds{\eta}\in
Q_{1,0}^{(i)}}\bigl\|\nabla
g_t(\bolds{\eta})\bigr\|+\sup_{\bolds{\eta}\in
Q_{2,0}^{(i)}}\bigl\|\nabla
g_t(\bolds{\eta})\bigr\| \Bigr),\nonumber
\end{eqnarray}
%
noting that $\partial g_t(\bolds{\eta})/\partial\bar{p}=
\partial g_t(\bolds{\eta})/\partial d$ does not exist at any point in
$\bar{B}_{\delta_1}(\bolds{\eta}_{i})$ with $d=d_0-1/2$.
As will be seen later, (\ref{A7}) together with
%
%
\begin{equation}
\label{A8} \max_{2\le t\le n}\mathrm{E} \Bigl(\sup_{\bolds{\eta}\in
Q_{v,0}^{(i)}}\bigl\|
\nabla g_t(\bolds{\eta})\bigr\| \Bigr) =O \bigl((\log
n)^{5/2} \bigr),\qquad v=1,2,
\end{equation}
constitutes a key step in the proof of (\ref{A9}).

To prove (\ref{A8}) for $v=1$, define
$g_t^{(L)}(\bolds{\eta})=\sqrt{n}\{n(1-B)\}^{d-d_0}A_{1,\bolds
{\theta}_0}^{-1}(B) \times
A_{2,\bolds{\theta}_0}(B)A_{1,\bolds{\theta
}}(B)A_{2,\bolds{\theta}}^{-1}(B)\varepsilon
_t-n^{d-d_0+1/2}\varepsilon_t$.
Then, $g_t^{(L)}(\bolds{\eta})=g_t(\bolds{\eta})$ for
$\bolds{\eta}\in Q_1^{(i)}$.
In addition,
$\nabla g_t^{(L)}(\bolds{\eta})= ((\nabla g_t^{(L)}(\bolds
{\eta}))_i )_{1\leq i\leq\bar{p}}$
satisfies
\begin{eqnarray}
\label{A10}\quad
&&
\frac{ (\nabla
g_t^{(L)}(\bolds{\eta}))_k}{n^{1/2}}\nonumber\\[-8pt]\\[-8pt]
&&\qquad=\cases{\displaystyle -\bigl\{n(1-B)\bigr\}^{d-d_0}
\frac{A_{2,\bolds{\theta}_0}(B)}{A_{1,\bolds{\theta
}_0}(B)A_{2,\bolds{\theta}}(B)} \varepsilon_{t-k}, \vspace*{2pt}\cr
\hspace*{126pt}\mbox{if $1\le k\le
p_1$},
\vspace*{2pt}\cr
\displaystyle \bigl\{n(1-B)\bigr\}^{d-d_0} \frac{A_{2,\bolds{\theta
}_0}(B)}{A_{1,\bolds{\theta}_0}(B)}
\frac{A_{1,\bolds{\theta}}(B)}{A^2_{2,\bolds{\theta}}(B)}
\varepsilon_{t+p_1-k}, \vspace*{2pt}\cr
\hspace*{126pt}\mbox{if $p_1 + 1 \le k
\le p_1 + p_2$},
\vspace*{2pt}\cr
\displaystyle \bigl\{n(1-B)\bigr\}^{d-d_0}
\bigl(\log n+\log(1-B)\bigr) \frac{A_{2,\bolds{\theta
}_0}(B)}{A_{1,\bolds{\theta}_0}(B)} \frac{A_{1,\bolds{\theta
}}(B)}{A_{2,\bolds{\theta}}(B)}
\varepsilon_{t}
\vspace*{2pt}\cr
\displaystyle\qquad{} -(\log n) n^{d-d_0}\varepsilon_t,
\qquad \mbox{if $k=\bar{p}$}.}\nonumber
\end{eqnarray}
Denote $\bolds{\eta}_i$ by $(\bolds{\theta}_i^{\mathsf
{T}},d_i)^{\mathsf{T}}$.
We first consider the case of $d_i\ge d_0-1/2$.
Write $(\nabla g_t^{(L)}(\bolds{\eta}))_k=\sum
_{s=1}^{t-1}c_{s,k}^{(n)}(\bolds{\eta})\varepsilon_{t-s}$,
and define $\bolds{\eta}_i^{(L)}=(\eta_{i,1}^{(L)},\ldots,\eta
_{i,\bar{p}}^{(L)})^{\mathsf{T}}:=(\bolds{\theta}_i^{\mathsf
{T}},\break d_0-1/2)^{\mathsf{T}}$.
Clearly, $c_{s,k}^{(n)}(\bolds{\eta})$ has continuous partial
derivatives, $D_{\mathbf j} c_{s,k}^{(n)}(\bolds{\eta})$ on $Q_1^{(i)}$.
Moreover, since for any $\bolds{\eta}=(\eta_1,\ldots,\eta_{\bar
{p}})^{\mathsf{T}}\in Q_1^{(i)}$, the
hypercube formed by $\bolds{\eta}^{(L)}_i$
and $\bolds{\eta}$ is included in $Q_1^{(i)}$, it follows from
(3.10) of \citet{Lai94} 
and the Cauchy--Schwarz inequality that for any $\bolds{\eta}\in
Q_1^{(i)}$, $1\le k\le\bar{p}$
and $t \geq2$,
%
%
\begin{eqnarray}
\label{A11} && \bigl\{\bigl(\nabla g_t^{(L)}(
\bolds{\eta})\bigr)_k-\bigl(\nabla g_t^{(L)}
\bigl(\bolds{\eta}_i^{(L)}\bigr)\bigr)_k
\bigr\}^2
\nonumber
\\
&&\qquad= \Biggl\{ \sum_{m=1}^{\bar{p}}\sum
_{{\mathbf j}\in J(m,\bar{p})} \int_{\eta_{i,j_m}^{(L)}}^{\eta
_{j_m}}\cdots
\int_{\eta_{i,j_1}^{(L)}}^{\eta_{j_1}} \sum_{s=1}^{t-1}
R_{s,{\mathbf j}}^{(k)}\varepsilon_{t-s} \,d\xi_{j_1}
\cdots d\xi_{j_m} \Biggr\}^2
\nonumber\\
&&\qquad\le2^{\bar{p}} \sum_{m=1}^{\bar{p}}\sum
_{{\mathbf j}\in J(m,\bar{p})} \Biggl\{ \int_{\eta_{i,j_m}^{(L)}}^{\eta_{j_m}}
\cdots\int_{\eta_{i,j_1}^{(L)}}^{\eta_{j_1}} \sum
_{s=1}^{t-1} R_{s,{\mathbf j}}^{(k)}
\varepsilon_{t-s} \,d\xi_{j_1} \cdots d\xi_{j_m}
\Biggr\}^2
\\[-4pt]
&&\qquad\le2^{\bar{p}} \sum_{m=1}^{\bar{p}}\sum
_{{\mathbf j}\in J(m,\bar{p})} \operatorname{vol}\bigl(Q_1^{(i)}(m,{
\mathbf j})\bigr) \nonumber\\[-4pt]
&&\hspace*{64.5pt}\qquad\quad{}\times\int\cdots\int_{Q_1^{(i)}(m,{\mathbf j})} \Biggl(\sum
_{s=1}^{t-1} R_{s,{\mathbf j}}^{(k)}
\varepsilon_{t-s} \Biggr)^2 \,d\xi_{j_1} \cdots d
\xi_{j_m},
\nonumber
\end{eqnarray}
where
\[
R_{s,{\mathbf j}}^{(k)} = R_{s,{\mathbf j}}^{(k)}(\xi
_{j_1},\ldots,\xi
_{j_m}) = D_{\mathbf j} c_{s,k}^{(n)}(\xi_1,\ldots,\xi_{\bar{p}})|_{\xi
_j=\eta_{i,j}^{(L)},j\notin{\mathbf j}}\
\]
and
\begin{eqnarray}
\hspace*{-4pt}&&
Q_1^{(i)}(m, {\mathbf j})\nonumber\\[-2pt]
\hspace*{-4pt}&&\qquad=\bigl\{(\eta_{j_1},\ldots,\eta_{j_m})\dvtx
\bigl(\eta_{i,1}^{(L)},\ldots,
\eta_{i,j_1-1}^{(L)},\eta_{j_1},\eta_{i,j_1+1}^{(L)},\ldots,\eta
_{i,j_2-1}^{(L)},
\nonumber\\[-2pt]
\hspace*{-4pt}&&\qquad\hspace*{86.4pt}\eta_{j_2},\eta_{i,j_2+1}^{(L)},\ldots,
\eta_{i,j_m-1}^{(L)}, \eta_{j_m},\eta_{i,j_m+1}^{(L)},\ldots,\eta
^{(L)}_{i,\bar{p}}\bigr)\in Q_1^{(i)}
\bigr\},\nonumber\\[-2pt]
\hspace*{-4pt}&&\eqntext{{\mathbf j}\in J(m, \bar{p}),}
\end{eqnarray}
is an $m$-dimensional partial sphere. Now, by (\ref{A2}), (\ref
{A10}), (\ref{A11})
and a change of the order of integration, we obtain for $1\le k\le\bar{p}$,
%
%
\begin{eqnarray}
\label{A12}
&&
\mathrm{E} \Bigl[ \sup_{\bolds{\eta}\in Q_1^{(i)}}
\bigl\{\bigl(\nabla
g_t^{(L)}(\bolds{\eta})\bigr)_k-\bigl(
\nabla g_t^{(L)}\bigl(\bolds{\eta}^{(L)}_i
\bigr)\bigr)_k\bigr\}^2 \Bigr] \nonumber\\[-4pt]
&&\qquad\le C \sum
_{m=1}^{\bar{p}}\sum_{{\mathbf j}\in J(m,\bar{p})}
\operatorname{vol}^{2}\bigl(Q_1^{(i)}(m,{\mathbf j})\bigr)
\nonumber\\[-9.5pt]\\[-9.5pt]
&&\qquad\quad\hspace*{61pt}{} \times\Biggl\{ \sum_{s=1}^{t-1}\max
_{{\mathbf j}\in J(m,\bar{p}), 1\le m\le\bar{p}} \sup_{\bolds{\eta}\in
Q_1^{(i)}} \bigl(D_{\mathbf j}c_{s,k}^{(n)}(
\bolds{\eta})\bigr)^2 \Biggr\}\nonumber\\[-4pt]
&&\qquad=O \bigl((\log n)^{5}
\bigr),
\nonumber
\end{eqnarray}
where the last relation follows from for any $1\le s\le t-1$ with $2\le
t\le n$,
$\sup_{\bolds{\eta}\in Q_1^{(i)}}|D_{\mathbf
j}c_{s,k}^{(n)}(\bolds{\eta})|\le
Cs^{-1/2}$ if $1\le k<\bar{p}$ and $\bar{p}\notin{\mathbf j}$;
$C (\log n)s^{-1/2}$ if $k=\bar{p}$ and $\bar{p}\notin{\mathbf j}$
or if $1\le k<\bar{p}$ and $\bar{p}\in{\mathbf j}$; and
$C (\log n)^2s^{-1/2}$ if $k=\bar{p}$ and $\bar{p}\in{\mathbf
j}$.\vspace*{1pt}
Similarly, for all $1\le t\le n-1$,
$\mathrm{E}((\nabla g_t^{(L)}(\bolds{\eta}_i^{(L)}))_k^2)\le
C \log n$ if $1\le k<\bar{p}$, and $C (\log n)^3$ if $k=\bar{p}$.
Combining these, with (\ref{A12}), yields that for $d_i\ge d_0-1/2$,
%
%
\begin{eqnarray}
\label{A15} \max_{2\le t\le n}\mathrm{E} \Bigl(\sup
_{\bolds{\eta}\in Q_{1,0}^{(i)}}\bigl\|\nabla g_t(\bolds{
\eta})\bigr\|^2 \Bigr) &\le&\max_{2\le t\le n}\mathrm{E} \Bigl(\sup
_{\bolds{\eta}\in Q_{1}^{(i)}}\bigl\|\nabla g_t^{(L)}(\bolds{
\eta})\bigr\|^2 \Bigr) \nonumber\\[-9pt]\\[-9pt]
&=& O \bigl((\log n)^5 \bigr).\nonumber\vadjust{\goodbreak}
\end{eqnarray}
Replacing $\bolds{\eta}^{(L)}_i$ by $\bolds{\eta}_i$ and
using an argument similar to that used in proving (\ref{A15}),
it can be shown that (\ref{A15}) is still valid in the case of $d_i<d_0-1/2$.
As a result, (\ref{A8}) holds for $v=1$.

Define
\[
g_t^{(R)}(\bolds{\eta})=(1-B)^{d-d_0}A_{2,\bolds
{\theta}_0}(B)A_{1,\bolds{\theta}_0}^{-1}(B)
A_{1,\bolds{\theta}}(B)A_{2,\bolds{\theta}}^{-1}(B)
\varepsilon_t-\varepsilon_t.
\]
Then, $g_t^{(R)}(\bolds{\eta})=g_t(\bolds{\eta})$ on
$\bolds{\eta}\in Q_2^{(i)}$.
Moreover, note that
%
%
\begin{equation}
\label{A16}\quad \bigl(\nabla g_t^{(R)}(\bolds{\eta})
\bigr)_k = \cases{\displaystyle  -(1-B)^{d-d_0}\frac{A_{2,\bolds{\theta
}_0}(B)}{A_{1,\bolds{\theta}_0}(B)A_{2,\bolds{\theta}}(B)}
\varepsilon_{t-k}, \vspace*{2pt}\cr
\qquad\mbox{if $1\le k\le p_1$},
\vspace*{2pt}\cr
\displaystyle (1-B)^{d-d_0}\frac{A_{2,\bolds{\theta}_0}(B)}{A_{1,\bolds
{\theta}_0}(B)} \frac{A_{1,\bolds{\theta}}(B)}{A^2_{2,\bolds
{\theta}}(B)}\varepsilon_{t+p_1-k},\vspace*{2pt}\cr
\qquad \mbox{if $p_1+1\le k\le p_1+p_2$},
\vspace*{2pt}\cr
\displaystyle (1-B)^{d-d_0}\bigl(\log(1-B)\bigr)\frac{A_{2,\bolds{\theta
}_0}(B)}{A_{1,\bolds{\theta}_0}(B)} \frac{A_{1,\bolds{\theta
}}(B)}{A_{2,\bolds{\theta}}(B)}
\varepsilon_t, \vspace*{2pt}\cr
\qquad \mbox{if $k=\bar{p}$}.}
\end{equation}
Using (\ref{A16}) in place of (\ref{A10}), we can prove (\ref{A8})
with $v=2$ in the same way as $v=1$.
The details are omitted to save space. Equipped with (\ref{A7}), (\ref
{A8}) and Chebyshev's inequality, we obtain
\begin{eqnarray*}
&&
\int_K^\infty P \bigl(R^c(\mu)
\bigr)\,d\mu\\
&&\qquad\le 2lq \int_K^\infty\max
_{1\le v\le2, 2\le t\le n} P \biggl(\sup_{\bolds{\eta}\in
Q_{v,0}^{(i)}}\bigl\|\nabla
g_t(\bolds{\eta})\bigr\| \ge\frac{\mu^{\ell_1/(2q)}}{2 \bar{p}^{1/2}}
\biggr) \,d\mu
\\
&&\qquad\le C(\log n)^{5/2}\int_K^\infty
\mu^{-\ell_1/(2q)} \,d\mu\\
&&\qquad\le C(\log n)^{5/2},
\end{eqnarray*}
where the last inequality is ensured by $\ell_1>2q$. Hence (\ref{A9})
is proved.

In view of (\ref{A6}) and (\ref{A9}), (\ref{A5}) follows by showing
that for all large $n$,
%
%
\begin{equation}
\label{A17} \mathrm{(I)}\le C.
\end{equation}
To establish (\ref{A17}), motivated by page 1543 of \citet{ChaIng11}, 
we are led to consider the $\bar{p}$-dimensional hypercube, $H^{(i)}$,
circumscribing
$\bar{B}_{\delta_1}(\bolds{\eta}_i)$.
Choose \mbox{$\mu\ge K$}. We then divide $H^{(i)}$ into sub-hypercubes,
$H^{(i)}_j(\mu)$
(indexed by $j$), of equal size,
each of which has an edge of length $2\delta_1(\lfloor\mu_1^{(\ell
_1+1)/(2q)}\rfloor+1)^{-1}$
and a circumscribed hypersphere of radius $\sqrt{\bar{p}}\delta
_1(\lfloor\mu_1^{(\ell_1+1)/(2q)}\rfloor+1)^{-1}$,
where $\lfloor a\rfloor$ denotes the largest integer $\le a$.
Let $G_j^{(i)}(\mu)=\bar{B}_{\delta_1}(\bolds{\eta}_{i})\cap
H_j^{(i)}(\mu)$ and $\{G_{v_j}^{(i)}(\mu), j=1,\ldots,m^*\}$ be the
collection of
nonempty $G_j^{(i)}(\mu)$'s. Then it follows that
%
%
\begin{equation}
\label{A18} \bar{B}_{\delta_1}(\bolds{\eta}_{i})=\bigcup
_{j=1}^{m^*} G_{v_j}^{(i)}(
\mu) \quad\mbox{and}\quad m^*\le C^*\mu^{{(\ell_1+1)\bar{p}}/({2q})},
\end{equation}
where $C^*>0$ is independent of $\mu$. In addition, we have
%
%
\begin{eqnarray}
\label{A19}
A(\mu):\!&=& \Biggl\{ \inf_{\bolds{\eta}\in\bar
{B}_{\delta_1}(\bolds{\eta}_i)} \sum
_{r=0}^{\ell q-1}g^2_{n\iota+1+rz_n}(
\bolds{\eta})<\mu^{-{1}/{q}}, R(\mu) \Biggr\}
\nonumber\hspace*{-35pt}
\\
&\subset&\bigcup_{j=1}^{m^*}\bigcap
_{r=0}^{\ell q-1} \Bigl\{ \inf
_{\bolds{\eta}\in G_{v_j}^{(i)}(\mu)}\bigl|g_{n\iota+1+rz_n}(\bolds
{\eta})\bigr|<
\mu^{-{1}/({2q})},
\nonumber\hspace*{-35pt}\\[-8pt]\\[-8pt]
&&\hspace*{41.6pt} \mathop{\sup_{\|\bolds{\eta}_a-\bolds{\eta}_b\|\le c_{\mu}}}_{
\bolds{\eta}_a,\bolds{\eta}_b\in\bar{B}_{\delta_1}(\bolds
{\eta}_i)} \bigl|g_{n\iota+1+rz_n}(\bolds{
\eta}_a)-g_{n\iota+1+rz_n}(\bolds{\eta}_b)\bigr|<2
\mu^{-{1}/({2q})} \Bigr\}\nonumber\hspace*{-35pt}\\
:\!&=& \bigcup_{j=1}^{m^*}
\bigcap_{r=0}^{\ell q-1} D_{j,r}(\mu).
\nonumber\hspace*{-35pt}
\end{eqnarray}
Let $\bolds{\eta}_a^{(j)}\in G_{v_j}^{(i)}(\mu)$,
$j=1,\ldots,m^*$, be arbitrarily chosen. Then for any
$\bolds{\eta}\in G_{v_j}^{(i)}(\mu)$, we have
\begin{eqnarray*}
\bigl|g_{n\iota+1+rz_n}\bigl(\bolds{\eta}_a^{(j)}
\bigr)\bigr| &\le&\bigl|g_{n\iota+1+rz_n}\bigl(\bolds{\eta}_a^{(j)}
\bigr)-g_{n\iota+1+rz_n}(\bolds{\eta})\bigr| +\bigl|g_{n\iota
+1+rz_n}(\bolds{
\eta})\bigr|
\\
&\le&\mathop{\sup_{\|\bolds{\eta}_a-\bolds{\eta}_b\|\le
c_{\mu}}}_{
\bolds{\eta}_a,\bolds{\eta}_b\in\bar{B}_{\delta_1}(\bolds
{\eta}_i)} \bigl|g_{n\iota+1+rz_n}(\bolds{
\eta}_a)\,{-}\,g_{n\iota+1+rz_n}(\bolds{\eta}_b)\bigr|\,{+}\,\bigl|g_{n\iota+1+rz_n}(
\bolds{\eta})\bigr|,
\end{eqnarray*}
where the last inequality is ensured by
$\|\bolds{\eta}_a^{(j)}-\bolds{\eta}\|\le
2\delta_1\sqrt{\bar{p}}(\lfloor\mu^{(\ell_1+1)/(2q)}\rfloor+1)^{-1}
\le c_{\mu}$, and hence,
$D_{j,r}(\mu)\subset S_{j,r}(\mu):=\{|g_{n\iota+1+rz_n}(\bolds{\eta
}_a^{(j)})|<3\mu^{-{1}/({2q})}\}$.
Combining this with (\ref{A19}) yields
%
%
\begin{equation}
\label{A21} \mathrm{(I)}\le\int_K^\infty\sum
_{j=1}^{m^*}P \Biggl(\bigcap
_{r=0}^{\ell q-1}S_{j,r}(\mu) \Biggr) \,d\mu.
\end{equation}

It is shown at the end of the proof that for some $s_0>0$ and all large $t$,
%
%
\begin{equation}
\label{A22} \inf_{\bolds{\eta}\in\bar{B}_{\delta_1}(\bolds
{\eta}_i)} \mathrm{E}^{{1}/{2}}\bigl\{\bigl(
\varepsilon_t(\bolds{\eta})-\varepsilon_t
\bigr)^2a_t^{-1}(d) t\bigr
\}>s_0.
\end{equation}
By letting
$K=\max\{1,[(6\sigma/(s_0\delta_0))\max\{(\ell q/(1-\iota
))^{d_0-1/2-L},1\}]^{2q}\}$,
(\ref{A17}) follows.
To see this, denote by ${\cal F}_t$
the $\sigma$-algebra generated by $\{\varepsilon_s, 1\le s\le t\}$, and
recall that\vadjust{\goodbreak}
$\varepsilon_t(\bolds{\eta})=\sum_{s=0}^{t-1}b_s(\bolds{\eta
})\varepsilon_{t-s}$,
with $b_0(\bolds{\eta})=1$.
We obtain, after some algebraic manipulations,
%
%
\begin{equation}
\label{A23} P \Biggl(\bigcap_{r=0}^{\ell q-1}S_{j,r}(
\mu) \Biggr)=\mathrm{E} \Biggl[\prod_{r=0}^{\ell q-2}I_{S_{j,r}(\mu)}P
\bigl(S_{j,\ell q-1}(\mu)\mid{\cal F}_{n+1-2z_n}\bigr) \Biggr]
\end{equation}
and
%
%
\begin{eqnarray}
\label{A24}\quad && P\bigl(S_{j,\ell q-1}(\mu)\mid{\cal F}_{n+1-2z_n}\bigr)
\nonumber\\[-2pt]
&&\qquad= P \biggl( M_1\bigl(\bolds{\eta}^{(j)}_a,
\mu\bigr)<\frac{\sum_{s=1}^{z_n-1}b_s(\bolds{\eta
}_a^{(j)})\varepsilon_{n+1-z_n-s}} {
\operatorname{var}^{1/2}((\varepsilon_{z_n}(\bolds{\eta
}_a^{(j)})-\varepsilon_{z_n})/\sigma)} <M_2\bigl(\bolds{
\eta}^{(j)}_a,\mu\bigr)\Big\mid\\[-2pt]
&&\hspace*{266.3pt}{\cal F}_{n+1-2z_n}
\biggr),
\nonumber\vspace*{-3pt}
\end{eqnarray}
where
\begin{eqnarray*}
M_i(\bolds{\eta},\mu)&=& \frac{(-1)^{i}3\mu^{-{1}/({2q})}a_n^{
1/2}(d)n^{-{1}/{2}}} {
\operatorname{var}^{1/2}((\varepsilon_{z_n}(\bolds{\eta
})-\varepsilon
_{z_n})/\sigma)}\\[-2pt]
&&{} -
\frac{\sum_{s=z_n}^{n-z_n}b_s(\bolds{\eta})\varepsilon
_{n+1-z_n-s}} {
\operatorname{var}^{1/2}((\varepsilon_{z_n}(\bolds{\eta
})-\varepsilon
_{z_n})/\sigma)},\qquad i=1,2.
\end{eqnarray*}
Since (\ref{A22}) yields that for $\mu\ge K$ and $n$ sufficiently large,
\begin{eqnarray*}
&& M_2\bigl(\bolds{\eta}_a^{(j)},\mu
\bigr)-M_1\bigl(\bolds{\eta}_a^{(j)},
\mu\bigr) \\[-2pt]
&&\qquad= 6\sigma\mu^{-{1}/({2q})}\mathrm{E}^{-1/2}\bigl\{\bigl(
\varepsilon_{z_n}\bigl(\bolds{\eta}_a^{(j)}
\bigr)-\varepsilon_{z_n}\bigr)^2a^{-1}_{z_n}(d)z_n
\bigr\} \frac{a_n^{1/2}(d)n^{-1/2}}{a_{z_n}^{
1/2}(d)z_n^{-1/2}}
\\[-2pt]
&&\qquad\le6\sigma\mu^{-{1}/({2q})}s_0^{-1}\max\biggl\{
\biggl(\frac{\ell q}{1-\iota} \biggr)^{d_0-1/2-L},1 \biggr\} \le
\delta_0
\end{eqnarray*}
and since
$\sum_{s=1}^{z_n-1}b_s^2(\bolds{\eta}_a^{(j)})=
\operatorname{var}((\varepsilon_{z_n}(\bolds{\eta}_a^{(j)})-\varepsilon
_{z_n})/\sigma)$,
it follows from (A1), (\ref{A23}) and (\ref{A24}) that
\begin{eqnarray*}
&&
P \Biggl(\bigcap_{t=0}^{\ell q-1}S_{j,r}(
\mu) \Biggr) \\[-2pt]
&&\qquad\le M_1 \biggl(6\sigma\mu^{-{1}/({2q})}s_0^{-1}
\max\biggl\{ \biggl(\frac{\ell q}{1-\iota} \biggr)^{d_0-1/2-L},1
\biggr\}
\biggr)^{\alpha_0} \mathrm{E} \Biggl(\prod_{r=0}^{\ell
q-2}I_{S_{j,r}(\mu)}
\Biggr).
\end{eqnarray*}
Moreover, since $\ell>q^{-1}(\iota^{-1}-1)$, we can repeat the same
argument $\ell q$ times to get
%
%
\begin{eqnarray}
\label{A25}
 P \Biggl(\bigcap_{r=0}^{\ell q-1}\!S_{j,r}(
\mu) \Biggr)
\,{\le}\, M_1^{\ell q} \biggl( \!6\sigma
\mu^{-{1}/({2q})}s_0^{-1}\max\biggl\{\! \biggl(
\frac{\ell q}{1-\iota} \biggr)^{\!d_0-1/2-L},1\! \biggr\} \!\biggr
)^{\!\alpha_0\ell q}.\hspace*{-45pt}\vadjust{\goodbreak}
\end{eqnarray}
By noting that the bound on the right-hand side of (\ref{A25}) is
independent of $j$
and $\ell>(2/\alpha_0)+\{(\ell_1+1)\bar{p}/(\alpha_0 q)\}$, it is
concluded from (\ref{A18}), (\ref{A21}) and (\ref{A25})
that for all large $n$, (I) is bounded above by
\begin{eqnarray*}
&&
C^*M_1^{\ell q}\bigl( 6\sigma s_0^{-1}
\max\bigl\{\bigl(\ell q/(1-\iota)\bigr)^{d_0-1/2-L},1 \bigr\}
\bigr)^{\alpha_0\ell q} \\
&&\qquad{}\times\int_K^{\infty}
\mu^{-({\alpha_0\ell}/{2}-{(\ell_1+1)\bar{p}}/({2q}))}\,d\mu\le C.
\end{eqnarray*}

To complete the proof of the lemma, it only remains to show (\ref
{A22}). Define
%
%
\begin{equation}
\label{A26} U_{t,\bolds{\theta}} = \frac{A_{2,\bolds{\theta
}_0}(B)}{A_{1,\bolds{\theta}_0}(B)} \frac{A_{1,\bolds{\theta
}}(B)}{A_{2,\bolds{\theta}}(B)}
\varepsilon_t = \sum_{s=0}^{t-1}\,d_s(
\bolds{\theta})\varepsilon_{t-s},
\end{equation}
noting that $d_0(\bolds{\theta})=1$. By (\ref{13})--(\ref{15}),
there exist $0<c_1,c_2<\infty$ such that
%
%
\begin{equation}
\label{A27} \sup_{\bolds{\theta}\in\Pi}\bigl|d_s(\bolds{
\theta})\bigr|\le c_1\exp(-c_2 s).
\end{equation}
%
Let $\Sigma_t(\bolds{\theta})=\mathrm{E}[(U_{1,\bolds{\theta
}},\ldots,U_{t,\bolds{\theta}})^{\mathsf{T}}(U_{1,\bolds
{\theta}},\ldots,U_{t,\bolds{\theta}})]$.
It can be shown that
there exists $C_0>0$ such that for all $t\ge1$,
%
%
\begin{equation}
\label{A31} \inf_{\bolds{\theta}\in\Pi}\lambda_{\min}\bigl(
\Sigma_t(\bolds{\theta})\bigr) \ge C_0
\sigma^{2}.
\end{equation}
%
Express $\varepsilon_t(\bolds{\eta})$ as
$(1-B)^{d-d_0}U_{t,\bolds{\theta}}=\sum
_{s=0}^{t-1}v_s(d)U_{t-s,\bolds{\theta}}$.
Then, there exists $G=G(L,U)>0$,
depending only on $L$ and $U$, such that for any $s\ge0$ and $d\in D$,
%
%
\begin{equation}
\label{A33} \bigl|v_s(d)\bigr|\le G(s+1)^{d_0-d-1}.
\end{equation}
Let $0<\eta^*<1/2$ be given. Straightforward calculations yield that
there exists
$C_{\eta^*}>0$ such that for any $s\ge0$ and $L\le d\le d_0-\eta^*$,
%
%
\begin{equation}
\label{A34} \bigl|v_s(d)\bigr|\ge C_{\eta^*}(s+1)^{d_0-d-1}.
\end{equation}
Let $\iota_1>0$ be small enough such that
%
%
\begin{eqnarray}
\label{A35} &0<\iota_1<\tfrac12-\eta^*,\qquad d_0-\tfrac12-
\iota_1>\underline{D}_i,&\nonumber\\[-8pt]\\[-8pt]
&d_0-\tfrac12+
\iota_1<\bar{D}_i
\quad\mbox{and}\quad C_0C^2_{\eta^*}\bigl(\log
\iota_1^{-1}\bigr)\iota_1^{2\iota_1}>2.&
\nonumber
\end{eqnarray}
Define
\begin{eqnarray*}
A_1&=&\bigl\{\bigl(\bolds{\theta}^{\mathsf{T}},d
\bigr)^{\mathsf{T}}\dvtx\bigl(\bolds{\theta}^{\mathsf{T}},d
\bigr)^{\mathsf{T}}\in\bar{B}_{\delta_1}(\bolds{\eta}_i),
\underline{D}_i\le d\le d_0-\tfrac12-
\iota_1 \bigr\},
\\
A_2&=&\bigl\{\bigl(\bolds{\theta}^{\mathsf{T}},d
\bigr)^{\mathsf{T}}\dvtx\bigl(\bolds{\theta}^{\mathsf{T}},d
\bigr)^{\mathsf{T}}\in\bar{B}_{\delta_1}(\bolds{\eta}_i),
d_0-\tfrac12-\iota_1< d< d_0-\tfrac12
\bigr\},
\\
A_3&=&\bigl\{\bigl(\bolds{\theta}^{\mathsf{T}},d
\bigr)^{\mathsf{T}}\dvtx\bigl(\bolds{\theta}^{\mathsf{T}},d
\bigr)^{\mathsf{T}}\in\bar{B}_{\delta_1}(\bolds{\eta}_i),
d= d_0-\tfrac12 \bigr\},
\\
A_4&=&\bigl\{\bigl(\bolds{\theta}^{\mathsf{T}},d
\bigr)^{\mathsf{T}}\dvtx\bigl(\bolds{\theta}^{\mathsf{T}},d
\bigr)^{\mathsf{T}}\in\bar{B}_{\delta_1}(\bolds{\eta}_i),
d_0-\tfrac12< d< d_0-\tfrac12+\iota_1
\bigr\},
\\
A_5&=&\bigl\{\bigl(\bolds{\theta}^{\mathsf{T}},d
\bigr)^{\mathsf{T}}\dvtx\bigl(\bolds{\theta}^{\mathsf{T}},d
\bigr)^{\mathsf{T}}\in\bar{B}_{\delta_1}(\bolds{\eta}_i),
d_0-\tfrac12+\iota_1\le d\le\bar{D}_i
\bigr\}.
\end{eqnarray*}
Then, (\ref{A22}) is ensured by
showing that there exist $\zeta_{i}>0, i=1,\ldots, 5$,
such that for all large $t$
%
%
\begin{equation}
\label{A36} \inf_{\bolds{\eta}\in A_i}\mathrm{E}\bigl\{\bigl(
\varepsilon_t(\bolds{\eta})-\varepsilon_t
\bigr)^2a_t^{-1}(d) t\bigr\} >
\zeta_{i},\qquad i=1,\ldots,5.
\end{equation}

To show (\ref{A36}) with $i=1$, we deduce from (\ref{A31}), (\ref
{A34}), (\ref{A35}) and a straightforward calculation that for all
$\bolds{\eta}\in A_1$
and $t\ge2$,
%
%
\begin{eqnarray}
\label{A37} \mathrm{E} \bigl(\varepsilon_t^2(
\bolds{\eta})a_t^{-1}(d)t \bigr) &\ge&
\frac{\sigma^2 C_0C_{\eta^*}^2}{t^{2(d_0-d)-1}}\sum
_{s=0}^{t-1}(s+1)^{2(d_0-d)-2}\nonumber\\[-8pt]\\[-8pt]
&\ge&\frac{\sigma^2C_0C_{\eta^*}^2\{1-(1/2)^{2\iota_1}\}
}{2(d_0-\underline{D}_i)-1}.\nonumber
\end{eqnarray}
In addition, it is clear that
$\sup_{\bolds{\eta}\in A_1}\mathrm{E}(\varepsilon_t^2a_t^{-1}(d) t)\le
\sigma^2/t^{2\iota_1}\rightarrow0$, as $t\rightarrow\infty$,
which, together with (\ref{A37}), yields (\ref{A36}) with $i=1$.
By (\ref{A35}), Taylor's theorem and an argument similar to that of
(\ref{A37}), we have for all $\eta\in A_2$ and sufficiently large $t$,
%
%
\begin{eqnarray}
\label{A38} \mathrm{E} \bigl(\varepsilon_t^2(
\bolds{\eta})a_t^{-1}(d) t \bigr) &\ge&
\frac{\sigma^2C_0C_{\eta^*}^2}{t^{2(d_0-d)-1}}\int_{\iota_1
t}^{t}x^{2(d_0-d)-2}\,dx
\nonumber\\
&\ge&\frac{\sigma^2C_0C_{\eta^*}^2}{2(d_0-d)-1}\bigl(1-\iota
_1^{2(d_0-d)-1}\bigr) \\
&\ge&
\sigma^2C_0C_{\eta^*}^2\bigl(\log
\iota_1^{-1}\bigr)\iota_1^{2\iota_1}
>2\sigma^2.
\nonumber
\end{eqnarray}
Moreover, $\sup_{\bolds{\eta}\in A_2}\mathrm{E}(\varepsilon_t^2
a_t^{-1}(d) t)\le\sigma^2$.
This and (\ref{A38}) together imply (\ref{A36}) for $i=2$.
Equation (\ref{A36}) for $i=3$ follows directly from (\ref{A31}) and
(\ref{A34}).
The details are skipped.
To show that (\ref{A36}) holds for $i=4$, we get from Taylor's theorem
and (\ref{A35}) that for all $\bolds{\eta}\in A_4$
and sufficiently large $t$,
\begin{eqnarray*}
\mathrm{E} \bigl(\varepsilon_t^2(\bolds{
\eta})a_t^{-1}(d) t \bigr)&=&\mathrm{E}\bigl(
\varepsilon_t^2(\bolds{\eta})\bigr) \ge
\sigma^2C_0C_{\eta^*}^2\int
_1^{\iota_1^{-1}}x^{2(d_0-d)-2} \,dx
\\
&\ge&\sigma^2C_0C_{\eta^*}^2\bigl(
\log\iota_1^{-1}\bigr)\iota_1^{2\iota_1}>2
\sigma^2.
\end{eqnarray*}
Hence, the desired conclusion follows.
For $\bolds{\eta}\in A_5$, define $W_t(\bolds{\eta})=\mathrm
{E}(\varepsilon_t(\bolds{\eta})-\varepsilon_t)^2$.
Then, it follows from (\ref{A27}), (\ref{A33}) and the Weierstrass
convergence theorem that
$W_t(\bolds{\eta})$ converges uniformly on $A_5$
to some function $W_{\infty}(\bolds{\eta})$.
Moreover, since $W_t(\bolds{\eta})$ is continuous on $A_5$, the
uniform convergence of $W_t(\bolds{\eta})$
ensures that $W_\infty(\bolds{\eta})$ is also continuous. In view
of (\ref{13})--(\ref{15}) and (\ref{A2}), we have
$(1-z)^{d-d_0}A_{2,\bolds{\theta}_0}(z)A^{-1}_{1,\bolds{\theta
}_0}(z)A_{1,\bolds{\theta}}(z)A^{-1}_{2,\bolds{\theta}}(z)\ne1$.
Hence, for each $\bolds{\eta}\in A_5$, \mbox{$W_t(\bolds{\eta})>0$}
for all sufficiently large $t$.
This, the continuity of $W_\infty(\bolds{\eta})$ and the
compactness of $A_5$ yield that there exists $\tilde{\theta}>0$
such that $\inf_{\bolds{\eta}\in A_5}W_\infty(\bolds{\eta})>
\tilde{\theta}$.\vadjust{\goodbreak}
By making use of this finding and the uniform convergence of
$W_t(\bolds{\eta})$ to $W_\infty(\bolds{\eta})$,
we obtain (\ref{A36}) with $i=5$.
This complete the proof of Lemma~\ref{lemma1}.
\end{pf*}

\begin{pf*}{Proof of Lemma~\ref{lemma2}}
Following the proof of Lemma~\ref{lemma1}, write
$\varepsilon_t(\bolds{\eta})-\varepsilon_t=\sum
_{s=1}^{t-1}b_s(\bolds{\eta})\varepsilon_{t-s}$.
Note first that $b_s(\bolds{\eta})$ has continuous partial
derivatives, $D_{\mathbf j}b_s(\bolds{\eta})$, on $\Pi\times D$.
By an argument similar to that used to
derive bounds for
${\sup_{\bolds{\eta}\in Q_1^{(i)}}}|D_{\mathbf
j}c_{s,k}^{(n)}(\bolds{\eta})|$
in Lemma~\ref{lemma1}, we have for all $s \geq1$,
%
%
\begin{equation}
\label{A39}\quad \sup_{\bolds{\eta}\in B_{j,v}}\bigl|D_{\mathbf j}b_s(
\bolds{\eta})\bigr|\le\cases{ Cs^{-1/2}, &\quad if $j=0$ and $\bar{p}
\notin{\mathbf j}$,
\vspace*{2pt}\cr
Cs^{-1/2}\log(s+1), &\quad if $j=0$ and $\bar{p}\in{\mathbf
j}$,
\vspace*{2pt}\cr
Cs^{vj-1/2}, &\quad if $j\ge1$ and $\bar{p}\notin{\mathbf j}$,
\vspace*{2pt}\cr
Cs^{vj-1/2}\log(s+1), &\quad if $j\ge1$ and $\bar{p}\in{\mathbf j}$.}
\end{equation}
Moreover, it follows from (\ref{A27}) and (\ref{A33}) that
%
%
\begin{equation}
\label{A40} \sup_{\bolds{\eta}\in B_{j,v}}\sum_{s=0}^{n-1}b_s^2(
\bolds{\eta})= \cases{ O(\log n), &\quad if $j=0$,
\vspace*{2pt}\cr
O\bigl(n^{2vj}
\bigr), &\quad if $j\ge1$.}
\end{equation}
In view of (\ref{A39}), (\ref{A40}), the compactness of $B_{j,v}$,
$j\ge0$ and (B.5) of \citet{ChaIng11}, 
we get
%
%
\begin{eqnarray}
\label{A41} && \mathrm{E} \Biggl\{\sup_{\bolds{\eta}\in B_{j,v}} \Biggl|\sum
_{t=1}^n\bigl(\varepsilon_t(
\bolds{\eta})-\varepsilon_t(\bolds{\eta}_0)
\bigr)\varepsilon_t \Biggr|^{q_1} \Biggr\}
\nonumber
\\
&&\qquad\le Cn^{{q_1}/{2}} \Biggl[ \Biggl\{\sup_{\bolds{\eta}\in
B_{j,v}}\sum
_{s=1}^{n-1}b^2_s(
\bolds{\eta}) \Biggr\}^{{q_1}/{2}}\nonumber\\[-8pt]\\[-8pt]
&&\qquad\quad\hspace*{36pt}{} + \Biggl\{\sum
_{s=1}^{n-1}\max_{{\mathbf j}\in J(m,\bar{p}),1\le m\le\bar{p}} \sup
_{\bolds{\eta}\in B_{j,v}}\bigl(D_{\mathbf j}b_s(\bolds{
\eta})\bigr)^2 \Biggr\}^{{q_1}/{2}} \Biggr]
\nonumber\\
&&\qquad= \cases{ O \bigl(n(\log n)^3 \bigr)^{{q_1}/{2}}, &\quad if $j=0$,
\cr
O \bigl(n^{1+2vj}(\log n)^2 \bigr)^{{q_1}/{2}}, &\quad if $j
\ge1$.}
\nonumber
\end{eqnarray}
By H\"{o}lder's inequality, the left-hand side of (\ref{23}) is
bounded above by
\begin{eqnarray*}
&& \Biggl\{\mathrm{E} \Biggl(\inf_{\bolds{\eta}\in B_{j,v}} a^{-1}_{n}(d)
\sum_{t=1}^n \bigl(\varepsilon_t(
\bolds{\eta})-\varepsilon_t(\bolds{\eta}_0)
\bigr)^{2} \Biggr)^{-{qq_1}/({q_1-q})} \Biggr\}^{({q_1-q})/{q_1}}
\\
&&\qquad{}\times\Biggl\{\mathrm{E} \Biggl(\sup_{\bolds{\eta}\in B_{j,v}}
\Biggl|\sum
_{t=1}^n \bigl(\varepsilon_t(
\bolds{\eta})-\varepsilon_t(\bolds{\eta}_0)
\bigr)\varepsilon_t \Biggr|^{q_1} \Biggr) \Biggr
\}^{{q}/{q_1}}\\
&&\qquad{}\times \bigl(n^{1+2(j-1)v}I_{\{j\ge1\}}+nI_{\{j=0\}}
\bigr)^{-q},
\end{eqnarray*}
which together with (\ref{A41}) and Lemma~\ref{lemma1}, gives the
desired conclusion.\vadjust{\goodbreak}~%
\end{pf*}

\begin{pf*}{Proof of Theorem~\ref{thm33}}
By a calculation similar but more complicated than that in the
proof of Theorem~\ref{thm31}, we obtain
%
%
\begin{eqnarray}
\label{A55}\qquad &&
y_{n+h}-\hat{y}_{n+h}(\hat{\bolds{
\eta}}_{n})-\sum_{s=0}^{h-1}
\underline{c}_{s}(\bolds{\eta}_{0})
\varepsilon_{n+h-s}\hspace*{-8pt}
\nonumber\\[-8pt]\\[-8pt]
&&\quad=\bigl(\nabla^{(1)} \varepsilon_{n+1}(\bolds{
\theta}_{0})\bigr)^{\mathsf{T}}\underline{L}_{h}(
\bolds{\eta}_{0}) (\hat{\bolds{\theta}}_{n}-
\bolds{\theta}_{0})
+\underline{\mathbf{c}}^{\mathsf{T}}_{h}
(\bolds{\eta}_{0})\mathbf{w}_{n+1, h}(
\hat{d}_{n}-d_{0})+r_{n},\hspace*{-8pt}
\nonumber
\end{eqnarray}
where $r_{n}$ satisfies $n \mathrm{E}(r^{2}_{n})=o(1)$, and
$((\nabla^{(1)} \varepsilon_{n+1}(\bolds{\theta}_{0}))^{\mathsf
{T}}\underline{L}_{h}(\bolds{\eta}_{0}),
\underline{\mathbf{c}}^{\mathsf{T}}_{h} (\bolds{\eta
}_{0})\mathbf{w}_{n+1, h})^{\mathsf{T}}$
and $n^{1/2}(\hat{\bolds{\eta}}_{n}-\bolds{\eta}_{0})$ are
asymptotically independent. The desired conclusion (\ref{318})
follows by a direct application of (\ref{A55}), (\ref{311})
and Theorem~\ref{thm21}.
\end{pf*}
\end{appendix}

\section*{Acknowledgments}

We would like to thank an Associate Editor and two anonymous referees
for their insightful and constructive comments, which greatly improve
the presentation of this paper.

\begin{supplement}
\stitle{Supplement to ``Moment bounds and mean squared prediction errors
of long-memory time series''}
\slink[doi]{10.1214/13-AOS1110SUPP} 
\sdatatype{.pdf}
\sfilename{aos1110\_supp.pdf}
\sdescription{The supplementary material contains a Monte Carlo
experiment of finite
sample performance of the CSS predictor and the proof of (\ref{29}).}
\end{supplement}


\printaddresses

\end{document}